\documentclass[sn-mathphys]{sn-jnl}% Math and Physical Sciences Reference Style
\jyear{2021}%

\theoremstyle{thmstyleone}%
\theoremstyle{thmstyletwo}%
\newtheorem{example}{Example}%
\newtheorem{remark}{Remark}%

\theoremstyle{thmstylethree}%

\begin{document}

\title[A spectral collocation method for elliptic PDEs]
{A spectral collocation method for elliptic PDEs in irregular domains with Fourier extension}


\author[1]{\fnm{Xianru} \sur{Chen}}\email{chenxianru@hust.edu.cn}
\equalcont{These authors contributed equally to this work.}

\author*[2]{\fnm{Li} \sur{Lin}}\email{linlimath@whut.edu.cn}
\equalcont{These authors contributed equally to this work.}

\affil[1]{\orgdiv{School of Mathematics and Statistics}, \orgname{Huazhong
University of Science and Technology}, \orgaddress{\city{Wuhan}, \postcode{430074}, \state{Hubei}, \country{P. R. China}}}

\affil*[2]{\orgdiv{College of Science}, \orgname{Wuhan University of Technology}, \orgaddress{\city{Wuhan}, \postcode{430070}, \state{Hubei}, \country{P. R. China}}}


\abstract{Based on the Fourier extension, we propose an oversampling collocation method for solving the elliptic partial differential equations with variable coefficients over arbitrary irregular domains. This method only uses the function values on the equispaced nodes, which has low computational cost and versatility. While a variety of numerical experiments are presented to demonstrate the effectiveness of this method, it shows that the approximation error fast reaches a plateau with increasing the degrees of freedom, due to the inherent ill-conditioned of frames.}




\maketitle

\section{Introduction}
The high-precision function approximation on multivariate domains is still a challenging work. One well-known method is the least squares approximation \cite{ZHOUTAO2020,Cohen2013}. The key of this method is to find a suitable famliy of orthogonal polynomials over the given domains. However, when dealing with irregular domains in higher dimensions or approximating singular functions, it is much harder, even impossible, to construct suitable orthogonal polynomials. We therefore suggest to use a slightly weaker concept, namely frames \cite{framebook2016,SIAMRev}. Frames are a more flexible tool that have not received much attention yet, especially in solving partial differential equations (PDEs). The aim of this paper is to use the Fourier frame to solve the following elliptic PDEs with variable coefficients:
\begin{equation}\label{PDE0}
\begin{cases}
 -\nabla \cdot (\alpha(\boldsymbol{x})\nabla U(\boldsymbol{x})) + \beta(\boldsymbol{x}) U(\boldsymbol{x}) = F(\boldsymbol{x}) \quad \mathrm{in} \quad \Omega,\\
 U(\boldsymbol{x}) = H(\boldsymbol{x}) \quad \mathrm{on} \quad \partial\Omega,\\
\end{cases}
\end{equation}
where $\Omega \subseteq \mathbb{R}^d$ is an irregular compact domain, it can be simple-connected or multi-connected, and the coefficient terms satisfy $\alpha(\boldsymbol{x}), \beta(\boldsymbol{x}) \in C(\Omega)$, $\alpha\geq\beta_0>0$, $\beta\geq0$.

There are two typical numerical methods for PDEs defined in irregular domains. The first one is the transformation method, which maps the irregular domain to the regular domain by explicit smooth mapping \cite{1980Orszag}. However, even simple constant coefficient PDEs general become variable coefficient PDEs after the mapping is applied. At the same time, the method is limited to the problems with smooth or fixed number of piecewise smooth boundaries.  The second one is the continuation method, namely the fictitious domain method, which embeds the irregular domain into a larger regular domain through a certain kind of extension \cite{2009LSH}. For instance, the zero extensions, the functions $\alpha$, $\beta$ and $F$ are simply set to zero in the extended domain. However, due to the low regularity of the extended problem, its approximation accuracy is limited to the first or second order. While for smooth extension, the functions $\alpha$, $\beta$ and $F$ are smoothly extended to a larger regular domain, and then a suitable variational formula for the extension problem is established. Thus, the extended solution is as smooth as the original solution \cite{1996Elghaoui,2020shenjieJSC,2021shenjieJSC}.

Under necessary assumptions, the smooth extension is available \cite{2002JCP}. Gu and Shen proposed a spectral Petrov-Galerkin method which encloses the irregular domain into a larger rectangular domain, namely rectangular embedding \cite{2020shenjieJSC}. One of the variational schemes is only suitable for Poisson equation, and the other one is suitable for general PDEs. Numerical experiments also show that the $L^2$ error can only reach around $\mathcal{O}(10^{-6})$ when the degree of freedom is about $100$.
And when the degree of freedom is quite small, the accuracy is divergent. Under the same degrees of freedom, the method we proposed in this paper can achieve higher convergence accuracy in the $L^\infty$-norm.
Further, Gu and Shen presented another spectral method, named circular embedding \cite{2021shenjieJSC}. The main advantage of this approach is that the extended two-dimensional problem can be decomposed into a sequence of one-dimensional differential equations by using polar transformation, but these systems are nested and cannot be solved in parallel. Meanwhile, the method is difficult to generalize to $d$ dimensions $(d>3)$ domain, where $\Omega$ need to be a simply connected smooth domain. However, the method we proposed in this paper can be generalized to high order and high dimensions naturally, since only the function point value information in $\Omega$ and on the boundary $\partial\Omega$ is used.

The Fourier extension (FE) method is closely related to fictitious domain methods for solving certain PDEs using Fourier basis \cite{Badea2003,Lyon2010HighorderUS}, the main difference being the approximation in the extension region. In the fictitious domain methods, the function is explicitly extended outside the domain of interest \cite{Penven2012OnTS,Badea2001Daripa,Astrakmantsev}. While in the FE technique, the approximation in the extension region is determined implicitly through solving a least squares problem. The convergence properties and numerical algorithms of the FE method have some mature results \cite{2014FCM,2011lyon,2016fast,2018fast,2020AZ}. In this paper, we propose a spectral collocation method for solving elliptic PDEs by using FE, and we present numerous numerical experiments. We can obtain spectral convergence by only using the function information in $\Omega$ and on $\partial\Omega$. Moreover, we observe that the error fast reaches a plateau with increasing freedom, particular for sufficiently smooth solutions.

The organization of this paper is as follows. In Section \ref{section2}, we briefly state the FE problem. In Section \ref{section3}, based on the FE, we develop a spectral collocation method for the second-order elliptic PDEs over arbitrary irregular domains. In Section \ref{section4}, we present some numerical experiments to demonstrate the effectiveness of this method, followed by some concluding remarks in Section \ref{section5}.

\section{$d$-dimensional Fourier extension}\label{section2}
In this section, we mainly state the core ideas of FE problem and describe the existing numerical methods briefly. We also give numerical experiments to assist some statements.

\subsection{Fourier frame and Fourier extension}
Let $\Omega \subseteq \mathbb{R}^d$ be an arbitrary domain which is compactly contained in a hypercube $R=[-T,T]^d$, $T>1$, $d\geq1$. Let $H$ is a separable Hilbert space over the field $\mathbb{C}$, and we denote $\langle \cdot, \cdot \rangle$ and $\|\cdot\|$ as the inner product and norm of $H$. Given a function $f$, the aim of FE is to find a Fourier series $\mathcal{F}$, which only uses the function value information of $f$ at equispaced nodes and makes $\|\mathcal{F}-f\|$ is minimized.

Let $\boldsymbol{l} = (l_1,l_2,...,l_d)$ be the $d$-dimensional integer index and let $I_\Lambda$ be the corresponding countable index set. Moreover, we assume that the degrees of freedom in each dimension are equal, that is, $-n\leq l_i\leq n$, $i=1,...,d$. Let $\boldsymbol{x} = (x_1,x_2,...,x_d)$, the tensor Fourier basis functions on $R$ are defined as
$P_\Lambda:=\{\phi_{\boldsymbol{l}}(\boldsymbol{x})\}_{\boldsymbol{l} \in I_\Lambda}= \{\exp(i\pi\boldsymbol{x}\cdot\boldsymbol{l}/ T)\}_{\mathbf{l} \in I_\Lambda}$,
and we denote $N_\Lambda:=\mid P_\Lambda\mid = N^d$, $N=2n+1$. Note that an orthonormal basis on $[-T,T]$ fails to constitute a basis when restricted to the smaller interval $[-1,1]$, it forms the so-called frame \cite{framebook2016,SIAMRev}. Hence the sequences
$\{\phi_{\boldsymbol{l}}(\boldsymbol{x})\}_{\boldsymbol{l} \in I_\Lambda}$ form a set of Fourier frames over $\Omega \subseteq R$. Moreover, we define the function space
$$\mathcal{G}_{N_\Lambda} = \textrm{span}\{\phi_{\boldsymbol{l}}(\boldsymbol{x})\}_{\boldsymbol{l} \in I_\Lambda}.$$
The FE problem is now formalized as finding an approximation
\begin{equation*}
\mathcal{F}_{N_\Lambda}(f) = \sum_{\boldsymbol{l} \in I_\Lambda } a_{\boldsymbol{l} }  \phi_{\boldsymbol{l}}(\boldsymbol{x})
\end{equation*}
such that
\begin{align}\label{FE}
\mathcal{F}_{N_\Lambda}(f) := \mathop{\min}_{\forall g\in \mathcal{G}_{N_\Lambda}}
\|{f(\boldsymbol{x})-g(\boldsymbol{x})} \|,\quad \boldsymbol{x} \in \Omega.
\end{align}
We refer to the $\mathcal{F}_{N_\Lambda}(f)$ as the Fourier extension of $f$, it is the orthogonal projection onto $\mathcal{G}_{N_\Lambda}$, and it is uniquely described by a set of coefficients
$\boldsymbol{a} \in \mathbb{C}^{N_\Lambda}$, which is the minimizer of an approximation algorithm
\begin{align}\label{coeff}
\boldsymbol{a} = \mathop{\arg\min}_{\forall \boldsymbol{c}\in\mathbb{C}^{N_\Lambda}}
\left\| f(\boldsymbol{x}) - \sum_{\boldsymbol{l} \in I_\Lambda } c_{\boldsymbol{l}}  \phi_{\boldsymbol{l}}(\boldsymbol{x})\right\|.
\end{align}

\begin{remark}
Note that the choice of frames general depends on the function being approximated. For smooth functions, we also can use the Chebyshev or Legendre frames \cite{2020AZ}. For the algebraic singular or logarithmic singular functions, we prefer the frame of polynomial plus modified polynomial, see \cite{SIAMRev}. This paper involves differential operation, so we choose Fourier frames for convenience.
\end{remark}

\subsection{Discrete Fourier extension}
In order to avoid complex integral operation in $\Omega$, we usually adopt the oversampling collocation method to solve \eqref{coeff}, and the premise of realizing the collocation method is to find a set of appropriate collocation nodes first. Let $M =\gamma N$ is even, we choose a set of equispaced nodes on $R$ with $M$ points per dimension, and this set is denoted as
\begin{equation*}
P_R = \left\{\left(\frac{2Tk_1}{M},\frac{2Tk_2}{M}, ... ,\frac{2Tk_d}{M}\right)
: k_i \in \mathbb{Z}, \quad -\frac{M}{2} \leq k_i \leq \frac{M}{2} , \quad i=1,...,d \right\},
\end{equation*}
then $ N_R:=\mid P_R\mid  = M^d$. Further, we restrict these nodes $P_R$ of $R$ to its subdomain $\Omega$, and we denote the set of nodes located in $\Omega$ as $P_\Omega$, that is, $P_\Omega = P_R \cap \Omega$, $N_\Omega:=\mid P_ \Omega\mid$. Here we always choose $\gamma > 1$ such that $N_\Omega> N_\Lambda$.

Assuming a linear indexing $\boldsymbol{x}_k$ of $P_\Omega$ from $1$ to $N_\Omega$ and $\phi_j(\boldsymbol{x})$ of $P_\Lambda$ from $1$ to $N_\Lambda$, the norm $||\cdot||$ is taken as a discrete summation over a set of collocation nodes on $\Omega$, then the minimization \eqref{FE} can be reformulated as a discrete least squares problem
\begin{equation}\label{disFE2}
\mathcal{\bar{F}}_{N_\Lambda}(f) := \mathop{\arg\min}_{\forall g\in \mathcal{G}_{N_\Lambda}}
\sum_{\boldsymbol{x} \in P_\Omega}(f(\boldsymbol{x})-g(\boldsymbol{x}))^2.
\end{equation}
We define $\mathcal{\bar{F}}_{N_\Lambda}(f)$ as the discrete Fourier extension of $f$. The discrete least squares problem \eqref{disFE2} can be written as the following linear system, i.e.,
\begin{equation}\label{matrix2}
A \boldsymbol{a}= \boldsymbol{b},\quad A \in \mathbb{C}^{N_\Omega \times N_\Lambda},
\quad \boldsymbol{b} \in \mathbb{C}^{N_\Omega},
\end{equation}
where
$$ A_{k,j} =\frac{1}{\sqrt{N_R}} \phi_j(\boldsymbol{x}_k), \quad \boldsymbol{b} _k = f(\boldsymbol{x}_k),
 \quad 1 \leq j \leq N_\Lambda,\quad 1\leq k \leq N_\Omega .$$
This is a full and exponentially ill-conditioned linear system, it can be regularized by using truncated singular values decomposition (tSVD) with a tolerance $\varepsilon>0$.

However, the tSVD method is computationally expensive, i.e., $\mathcal{O}(N_\Lambda^3)$. In order to overcome this problem, Matthysen and Huybrechs proposed a fast and robust algorithm for the computation of FE, namely AA algorithm \cite{2016fast,2018fast,PHDMatthysen}. They found that it is possible to filter out the part that makes the system ill-conditioned by multiplying a factor on both sides of the linear system, i.e., one can transform the original ill-conditioned system $A$ in \eqref{matrix2} into a well-posed low-rank system $A-A^3$. The rank of matrix $A-A^3$ is determined by the size of plunge region, that is, the number of singular values whose values fall between the interval $(\varepsilon, 1-\varepsilon)$. When $d=1$, the AA algorithm reduces the amount of operations from $\mathcal{O}(N^3)$ to $\mathcal{O}(N\log^2 N)$ \cite{2016fast}. When $d=2$, the whole calculation amount of AA algorithm is $\mathcal{O}(N^2_\Lambda\log^2 N_\Lambda)$ \cite{2018fast}. Moreover, due to the connection with the trigonometric polynomials, the results of the size of plunge region can only be generalized to the Chebyshev frames \cite{2020AZ}. However, we only know the results when $d = 1, 2$, there is no fast algorithms for more higher dimensional frame approximation ($d\geq3$). Further, the AA algorithm is a particular case of the AZ algorithm \cite{2020AZ}, and this method is applicable as long as the singular values profile of matrix $A$ shows an exponential decay trend.

\begin{remark}
For one-dimensional FE, one usually takes double oversampling, i.e., $N_\Omega / N_\Lambda=2$. For high-dimensional FE, it is difficult to guarantee that the oversampling ratio $N_\Omega / N_\Lambda$ is a fixed constant. In Table \ref{first-table}, we give the number of collocation nodes on a diamond domain $\Omega_D$, where its vertices are $(1, 0)$, $(0, 1)$, $(-1, 0)$, $(0,-1)$. It shows that the value of $N_\Omega / N_\Lambda$ fluctuates slightly and changes quite small as $N$ increases. Hence we ignore the specific influence of ratio $N_\Omega / N_\Lambda$ on the approximation accuracy, and we also fix $T=2$, $\varepsilon=10^{-14}$ in this paper.
\end{remark}

\begin{table}[h]
\begin{center}
\begin{minipage}{\textwidth}
\caption{The number of collocation nodes used to discretize domain $\Omega_D$}
\label{first-table}
\begin{tabular}{@{}llllllllll@{}}
\toprule
$N$ & 10 & 15 & 20 & 25 & 30 & 35 & 40 & 45 & 50 \\
\midrule
$N_\Omega$ & 180 & 420 & 760 & 1200 & 1740 & 2380 & 3120 & 3960 & 4900 \\
$ \frac{N_\Omega}{N_\Lambda}$ & 1.8000 & 1.8667 & 1.9000 & 1.9200 & 1.9333 & 1.9429 & 1.9500 & 1.9556 & 1.9600 \\
\botrule
\end{tabular}
\end{minipage}
\end{center}
\end{table}

\begin{figure}[htbp]
\centering
{\includegraphics[width=5.8cm]{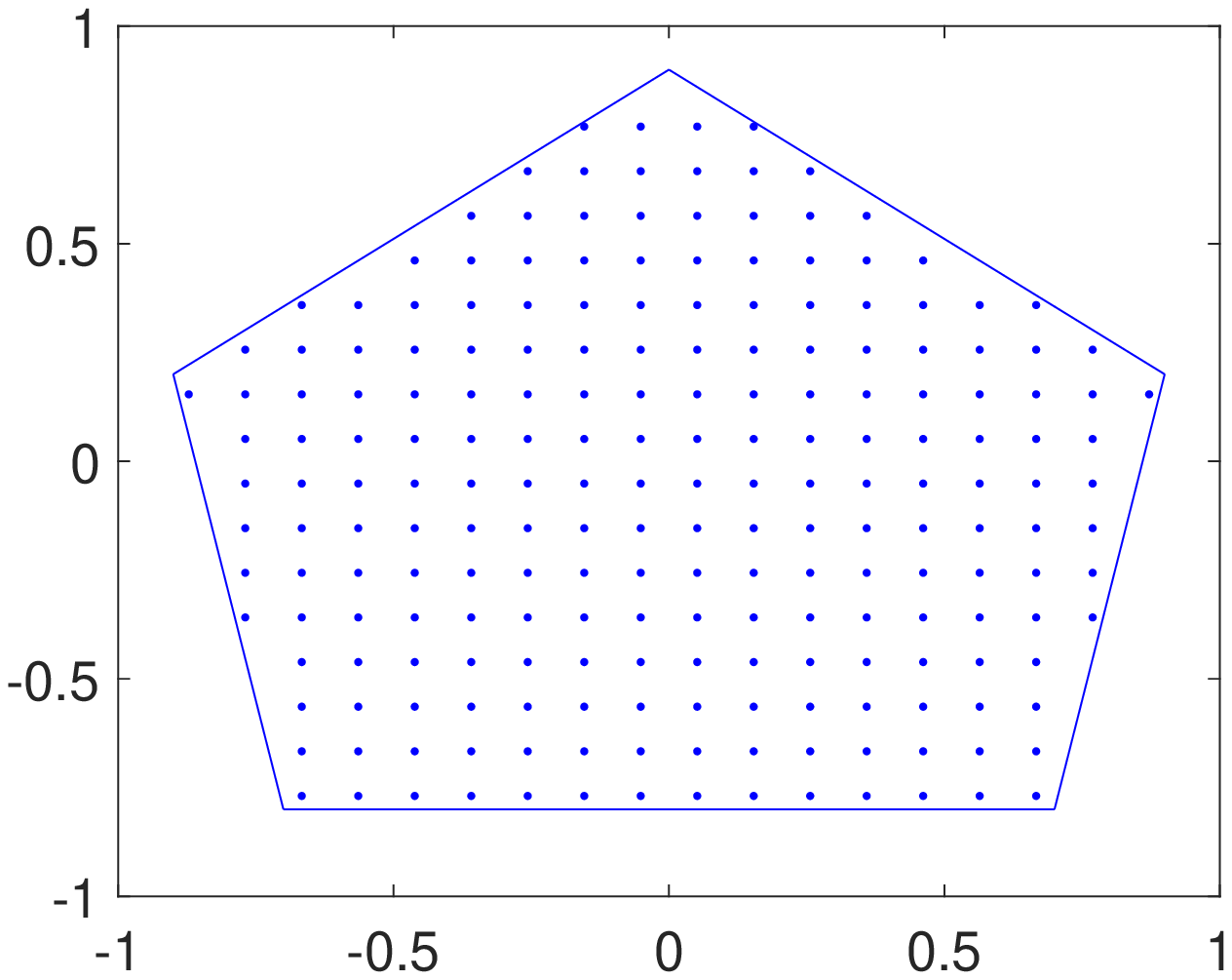}}
{\includegraphics[width=5.8cm]{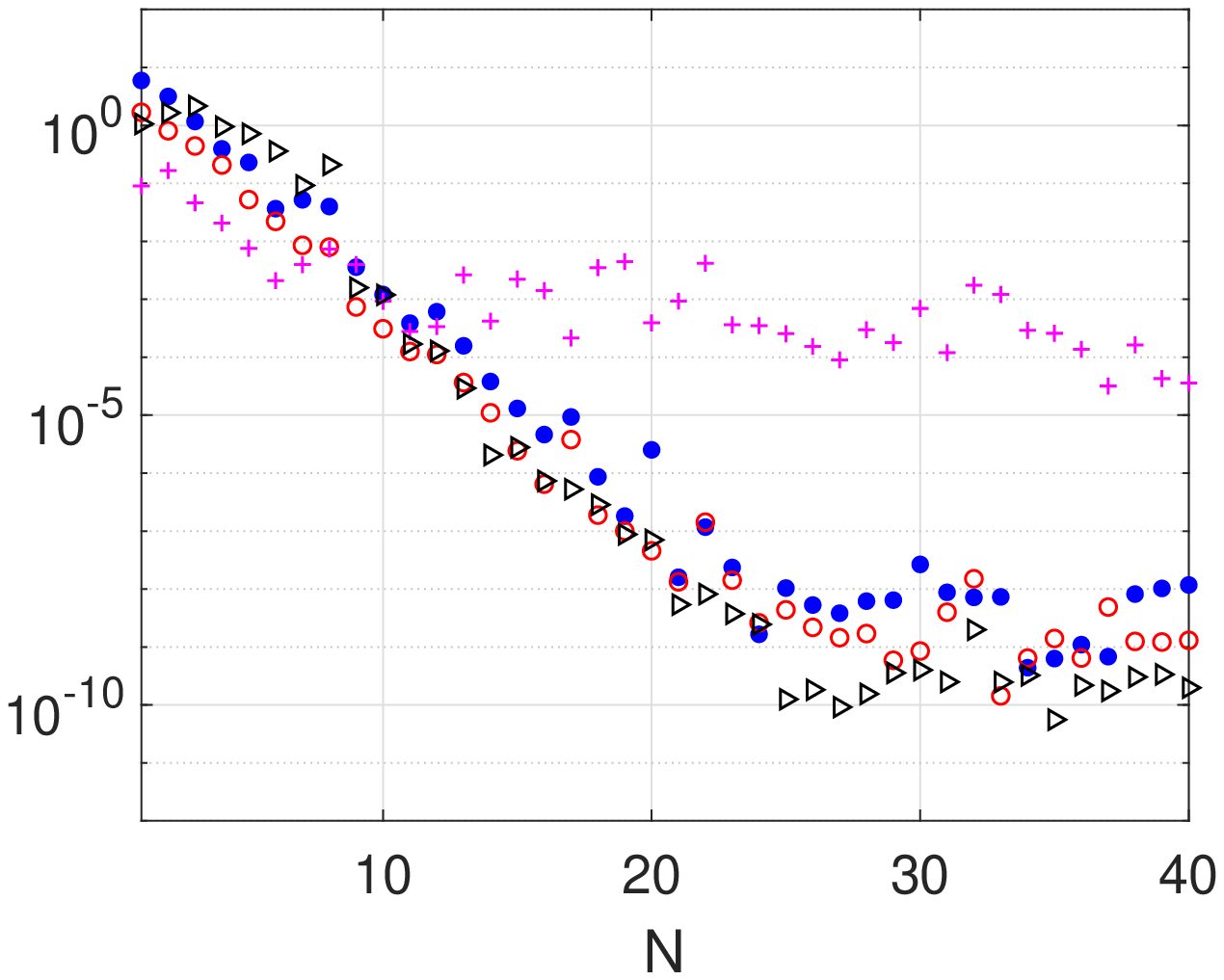}}
{\includegraphics[width=5.8cm]{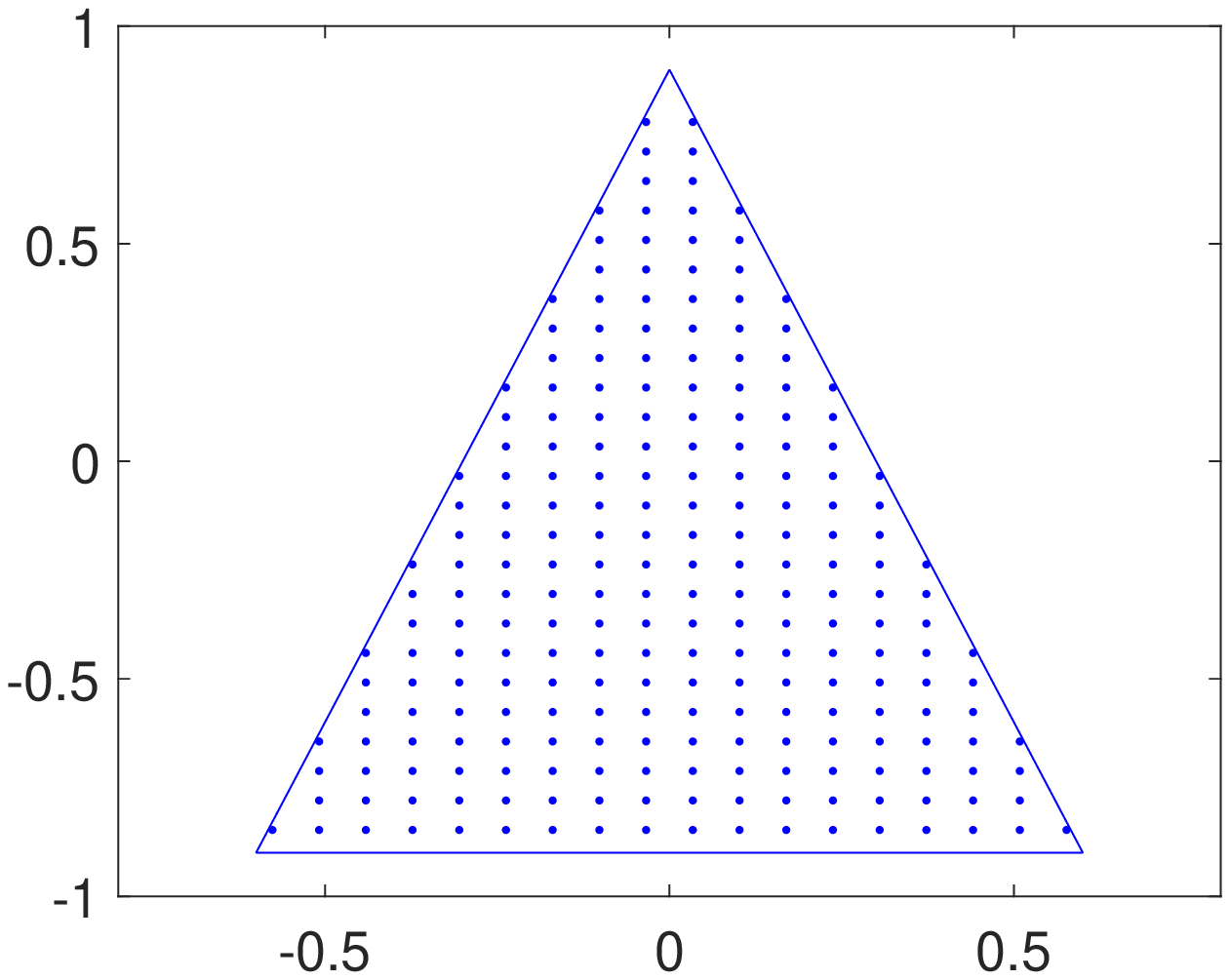}}
{\includegraphics[width=5.8cm]{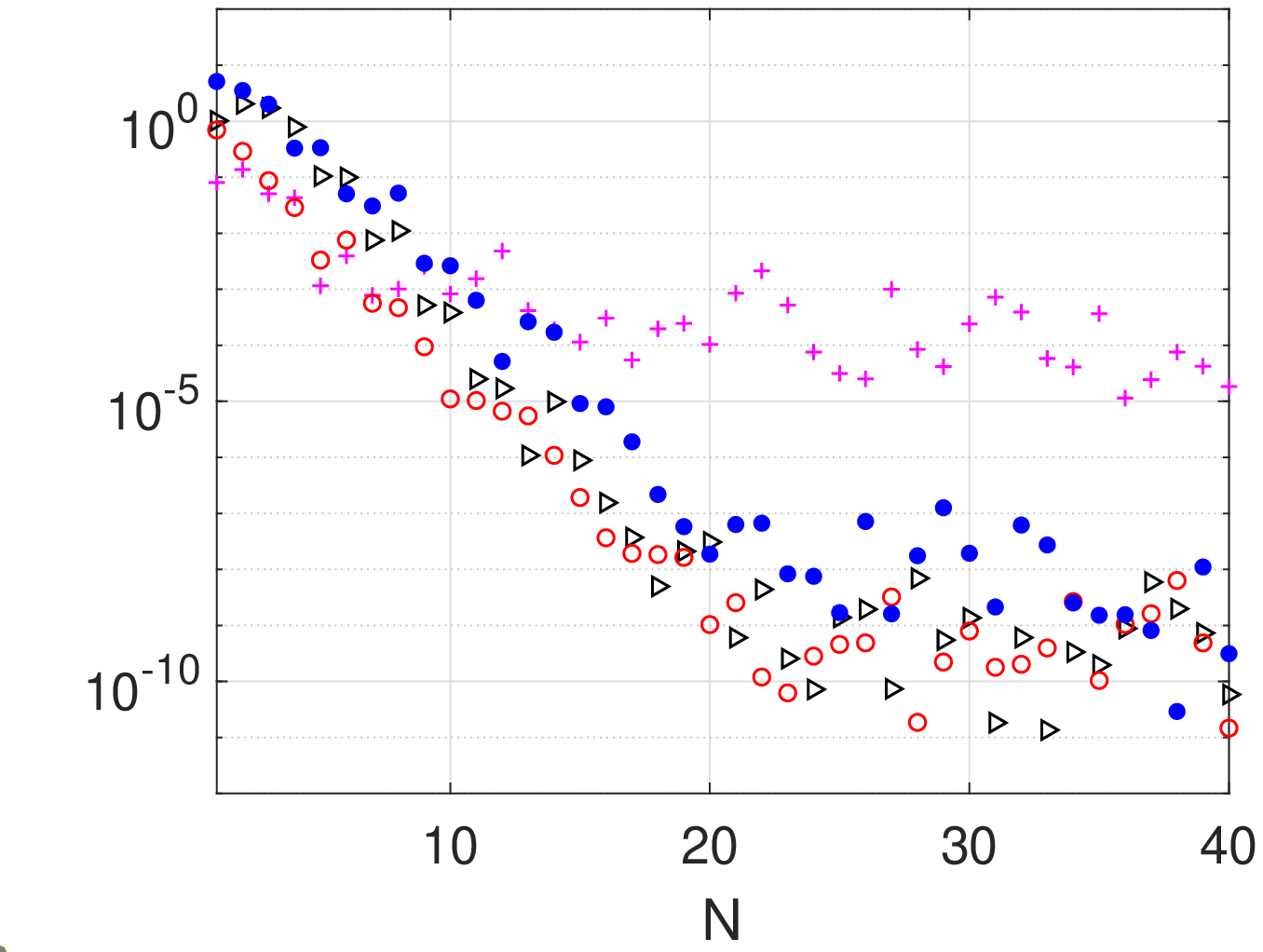}}
{\includegraphics[width=5.8cm]{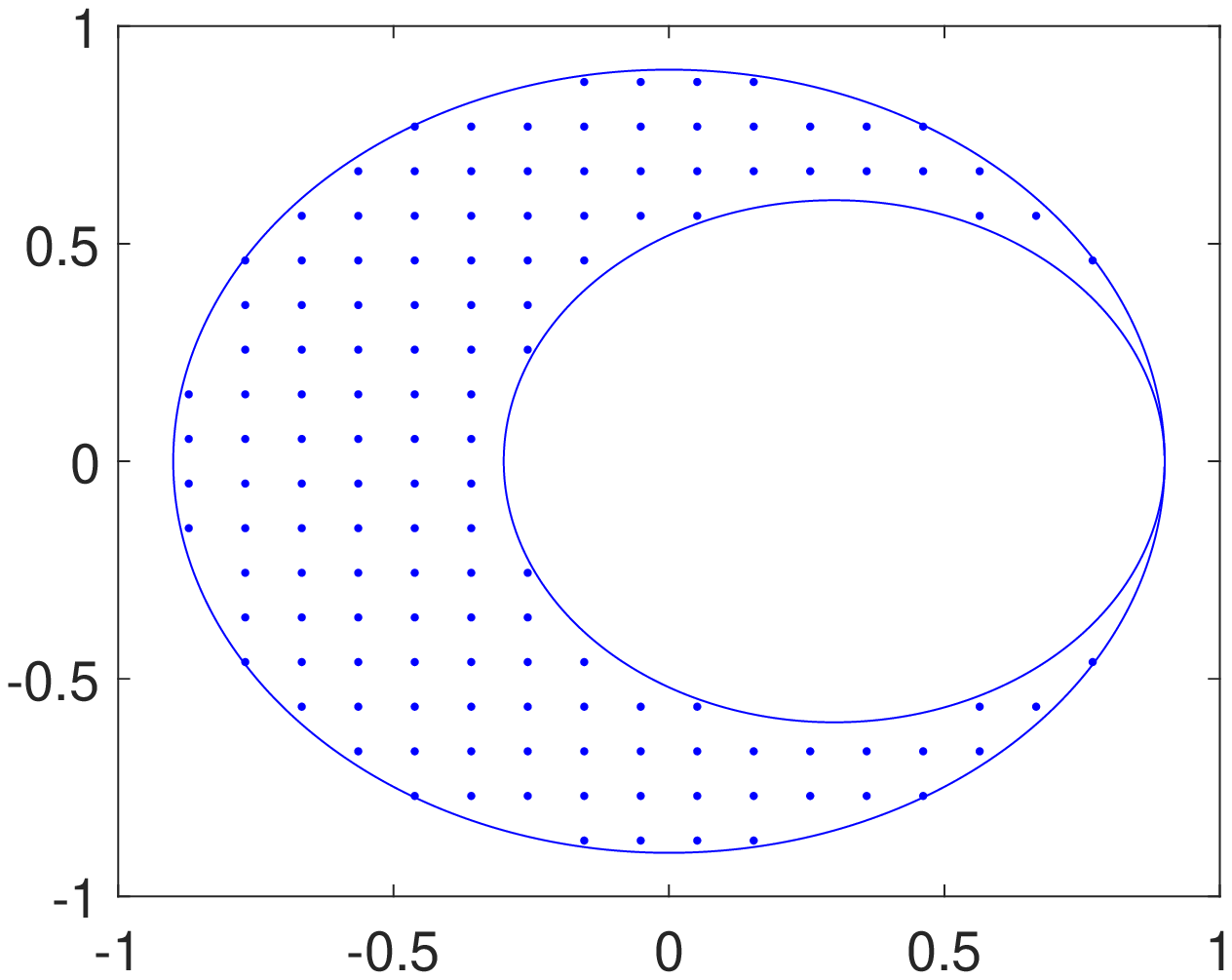}}
{\includegraphics[width=5.8cm]{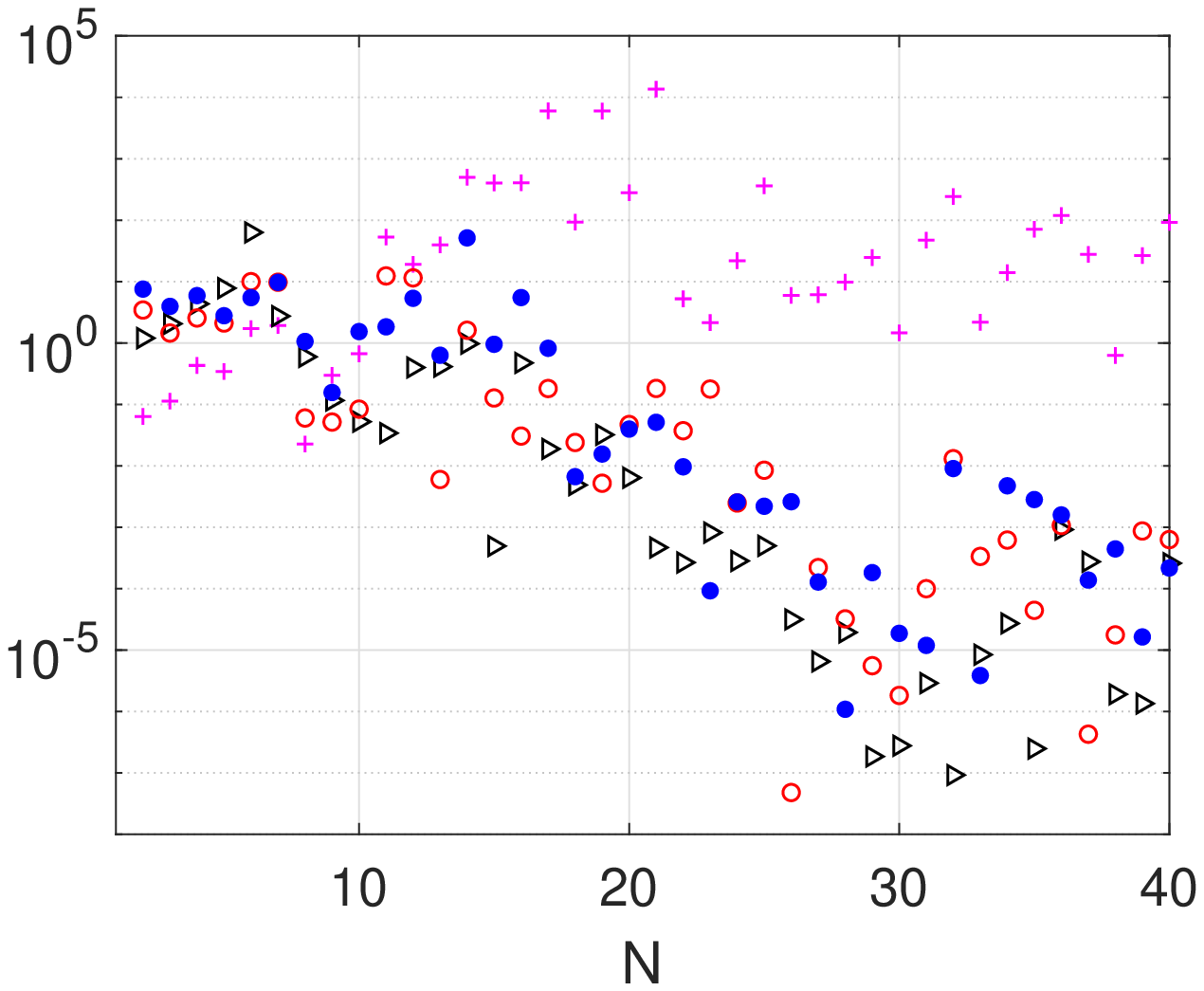}}
\caption{Left: The discretization nodes for domain $\Omega_P$, $\Omega_T$, $\Omega_L$ when $N=21$. Right: The log plot of the maximum error of functions $f_1$ (magenta pluses), $f_2$ (red circles), $f_3$ (black triangles), $f_4$ (blue dots) corresponding to the domains on the left side}
\label{2d-function-ex}
\end{figure}

In order to better understand the FE problem in two-dimensional case, we show the collocation nodes in three domains $\Omega_P$, $\Omega_T$, $\Omega_L$ and give the maximum error of four functions in Figure \ref{2d-function-ex}. The domain $\Omega_P$ is a pentagon with vertices $(0,0.9)$, $(-0.9,0.2)$, $(-0.7,-0.8)$, $(0.7,-0.8)$ and $(0.9,0.2)$. The domain $\Omega_T$ is a triangle with vertices $(0,0.9)$, $(-0.6,-0.9)$, $(0.6,-0.9)$. And the domain $\Omega_{L}=\{(x,y): x^2/0.9^2 + y^2/0.9^2-1 \leq 0,(x-0.3)^2/ 0.6^2 + y^2/0.6^2 -1 \geq 0 \}$. The functions are $f_1 =\mid xy\mid^3$, $f_2=1/((x-1.1)^2+(y-1.1)^2)^{3/2}$, $f_3=\cos(5x+y)\sin(x-3y)$ and $f_4= \exp(x+2y)$.
Due to the near-linear dependence of the truncated frames system, we observe that the error for analytic functions in $\Omega_P, \Omega_T$ will reach about $\mathcal{O}(10^{-10}) \sim \mathcal{O}(10^{-9})$ plateau, and there is no further improvement trend as $N$ increases. It should be noted that when using the Fourier frames to solve PDE problems, a similar phenomenon also occurs, although there is no specific theoretical analysis results. In particular, for domains with sharp corners, like $\Omega_L$, it can be seen that the approximation effects are generally poor. At this time, we need to add an appropriate number of collocation nodes at the sharp corners to improve the approximation accuracy \cite{PHDMatthysen}. The number of nodes should not be too large to affect the implementation of the AZ algorithm.

\unskip

These numerical experiments demonstrate that the approximation accuracy is affected by several factors, such as the regularity of functions, the shape of $\Omega$, the number of collocation nodes et al. It is difficult to obtain the optimal convergence accuracy and there is no convergence analysis. We need to maintain a balance between the amount of calculation and accuracy.

\section{The collocation method for elliptic PDEs}\label{section3}
Matthysen has promoted the AA algorithm that proposed in \cite{2016fast,2018fast,2020AZ}, so as to avoid the complexity of other domain-independent methods. This modification algorithm works based on the fact that when the collocation matrix is extended with some extra rows or columns, while satisfying the two requirements imposed in \cite{PHDMatthysen}, the singular value profile still holds. This makes the fast algorithm suitable to some extent for various problems that depend on function approximation, in particular includes the solution of elliptic boundary value problems with constant coefficient differential operators.
In this section, we mainly consider how to use the FE technique to numerically solve the variable coefficient elliptic PDEs \eqref{PDE0}.

\subsection{Discretization of the PDEs}
We use the Fourier frames $\{\phi_{\boldsymbol{l}} (\boldsymbol{x})\}$ defined on $\Omega$ to approximate the solutions of PDEs, i.e.,
\begin{equation*}
U(\boldsymbol{x}) \approx U_{N_\Lambda}(\boldsymbol{x})
:=\sum_{\boldsymbol{l} \in I_\Lambda } u_{\boldsymbol{l} }  \phi_{\boldsymbol{l}}(\boldsymbol{x}).
\end{equation*}
Then we employ the oversampling collocation method to discrete the PDEs and build the corresponding linear system. Let $\{\boldsymbol{x}_1, \boldsymbol{x}_2,...,\boldsymbol{x}_{N_I}\}$ be the nodes in the interior of $\Omega$ and let \{$\tilde{\boldsymbol{x}}_1, \tilde{\boldsymbol{x}}_2,...,\tilde{\boldsymbol{x}}_{N_B}\}$ be the nodes on the boundary $\partial\Omega$. Except for the requirement of oversampling, i.e., $N_I + N_B > N_\Lambda$, there is no clear requirement on the size of the number.
By imposing internal and boundary conditions at these nodes, we get
\begin{equation*}
\begin{split}
& \beta(\boldsymbol{x}_p) U_{N_\Lambda}(\boldsymbol{x}_p)
-\partial_x\alpha(\boldsymbol{x}_p)\partial_xU_{N_\Lambda}(\boldsymbol{x}_p)
-\partial_y\alpha(\boldsymbol{x}_p)\partial_yU_{N_\Lambda}(\boldsymbol{x}_p)\\
& \qquad - \alpha(\boldsymbol{x}_p) (\partial_{xx}U_{N_\Lambda}(\boldsymbol{x}_p)
 + \partial_{yy}U_{N_\Lambda}(\boldsymbol{x}_p))
= F(\boldsymbol{x}_p), \quad 1\leq p \leq N_I,\\
\end{split}
\end{equation*}
and
\begin{equation*}
U_{N_\Lambda}(\tilde{\boldsymbol{x}}_q) = H(\tilde{\boldsymbol{x}}_q),  \quad 1\leq q \leq N_B.
\end{equation*}
Let
\begin{equation*}
\begin{split}
&\boldsymbol{f}=(F(\boldsymbol{x}_1),...,F(\boldsymbol{x}_{N_I}))^T \in \mathbb{C}^{N_I},
\quad
\boldsymbol{h}=(H(\tilde{\boldsymbol{x}}_1),...,H(\tilde{\boldsymbol{x}}_{N_B}))^T \in \mathbb{C}^{N_B}, \\
&\boldsymbol{\alpha_1}= (\partial_x\alpha(\boldsymbol{x}_1),...,\partial_x\alpha(\boldsymbol{x}_{N_I}))^T \in \mathbb{C}^{N_I},
\quad
\boldsymbol{\alpha_2}=  (\partial_y\alpha(\boldsymbol{x}_1),...,\partial_y\alpha(\boldsymbol{x}_{N_I}))^T \in \mathbb{C}^{N_I}, \\
&\boldsymbol{\alpha_3}= (\alpha(\boldsymbol{x}_1),...,\alpha(\boldsymbol{x}_{N_I}))^T \in \mathbb{C}^{N_I},
\quad
\boldsymbol{\beta}=(\beta(\boldsymbol{x}_1),...,\beta(\boldsymbol{x}_{N_I}))^T \in \mathbb{C}^{N_I}.
\end{split}
\end{equation*}
Let $A1, A2, A3, A4, A5 \in \mathbb{C}^{N_I\times N_\Lambda}$ with entries
\begin{equation*}
\begin{split}
& (A1)_{j,k} = \phi_k(\boldsymbol{x}_j),
\quad  (A2)_{j,k} =\partial_x\phi_k(\boldsymbol{x}_j),
\quad (A3)_{j,k} =\partial_y\phi_k(\boldsymbol{x}_j),\\
& (A4)_{j,k} = \partial_{xx}\phi_k(\boldsymbol{x}_j),
\quad (A5)_{j,k} = \partial_{yy}\phi_k(\boldsymbol{x}_j),\quad 1\leq j \leq N_I,\quad 1\leq k\leq N_\Lambda,
\end{split}
\end{equation*}
and let $B \in \mathbb{C}^{N_B \times N_\Lambda}$ with entry
$$ B_{h,k} = \phi_k(\tilde{\boldsymbol{x}}_h), \quad 1\leq h \leq N_B,\quad 1\leq k\leq N_\Lambda.$$
Then we can establish the following rectangular linear system
\begin{equation}\label{PDE5}
\left(
  \begin{array}{ccc}
\boldsymbol{\beta}\circ A1 -\boldsymbol{\alpha_1} \circ A2 - \boldsymbol{\alpha_2} \circ A3 -\boldsymbol{\alpha_3} \circ (A4+A5) \\
B \\
  \end{array}
\right)
\boldsymbol{u}=
\left(
  \begin{array}{ccc}
   \boldsymbol{f} \\
    \boldsymbol{h} \\
  \end{array}
\right),
\end{equation}
where the circle $\circ$ represents the Hadamard product between the vector and the matrix.
We denote the coefficient matrix on the left side of \eqref{PDE5} as
$\boldsymbol{P} \in \mathbb{C}^{(N_I+N_B)\times N_\Lambda}$ and the point value vector on the right side as $\boldsymbol{F} \in \mathbb{C}^{N_I+N_B}$. Solving this full and ill-conditioned linear system
\begin{equation}\label{PDE6}
\boldsymbol{P} \boldsymbol{u} = \boldsymbol{F},
\end{equation}
we can obtain the numerical solution $U_N(\boldsymbol{x})$ of PDEs \eqref{PDE0}.

\subsection{Collocation nodes for solving PDEs}
Once the Fourier frames determined, the entries of the linear system \eqref{PDE5} only depend on the choice of collocation nodes, i.e., the shape of the domains plays an important role in the behavior of the solution. The hypercube $R$ is discretized with equispaced grids, and we restrict these grids to the interior of the irregular domain $\Omega$ to get the collocation nodes inside $\Omega$. The boundary is approximated with a set of discrete points lying on $\partial\Omega$. In practice, it is best for these points to be uniformly distributed across the boundary. In two dimensions, this can be accomplished easily by equally spacing points along an arc length parametrization of the curve. In three dimensions, however, equally distributing the points around a surface is more challenging.

We also need to concern the density of boundary nodes, i.e., the value of $N_B$. When using an insufficient number of nodes on the boundary, the accuracy suffers, while too many nodes can drive up the computational cost. Through a large number of numerical experiments, we find that it is generally enough to make $N_B= K N$, $K$ is generally  an integer greater than or equal to three. In fact, for single connected domains or multi-connected domains, the value $(N_I+N_B)/N_\Lambda$ gradually decreases as $N$ increases. Here we do not show the numerical experiments.

\section{Numerical  experiments}\label{section4}
To demonstrate the effectiveness of our method, we implement the algorithms in MATLAB and apply it to some examples already studied in the literature \cite{2009LSH,2020shenjieJSC}. Under the same degree of freedom, the proposed method has a more accurate solution. When discretizing a region, we use blue dots to represent the internal nodes and red dots to represent the boundary nodes.

\begin{example}\label{example1}
(Constant coefficient PDEs)~We set $\alpha(x,y) = 1$, $\beta(x,y) = 10$ in \eqref{PDE0}. Let the exact solution be $U(x,y)=\exp{\left(-(x^2+y^2)/2\right)}$, the pentagon doamin $\Omega_P$ is defined in Section \ref{section2}. Note that \eqref{PDE0} is a constant coefficient PDE, and then the linear system \eqref{PDE6} can be solved by the AA algorithm \cite{PHDMatthysen} as long as the number of boundary nodes satisfies $N_B = \mathcal{O}(\log N_\Lambda)$.

\begin{figure}[htbp]
\centering
{\includegraphics[width=5.8cm]{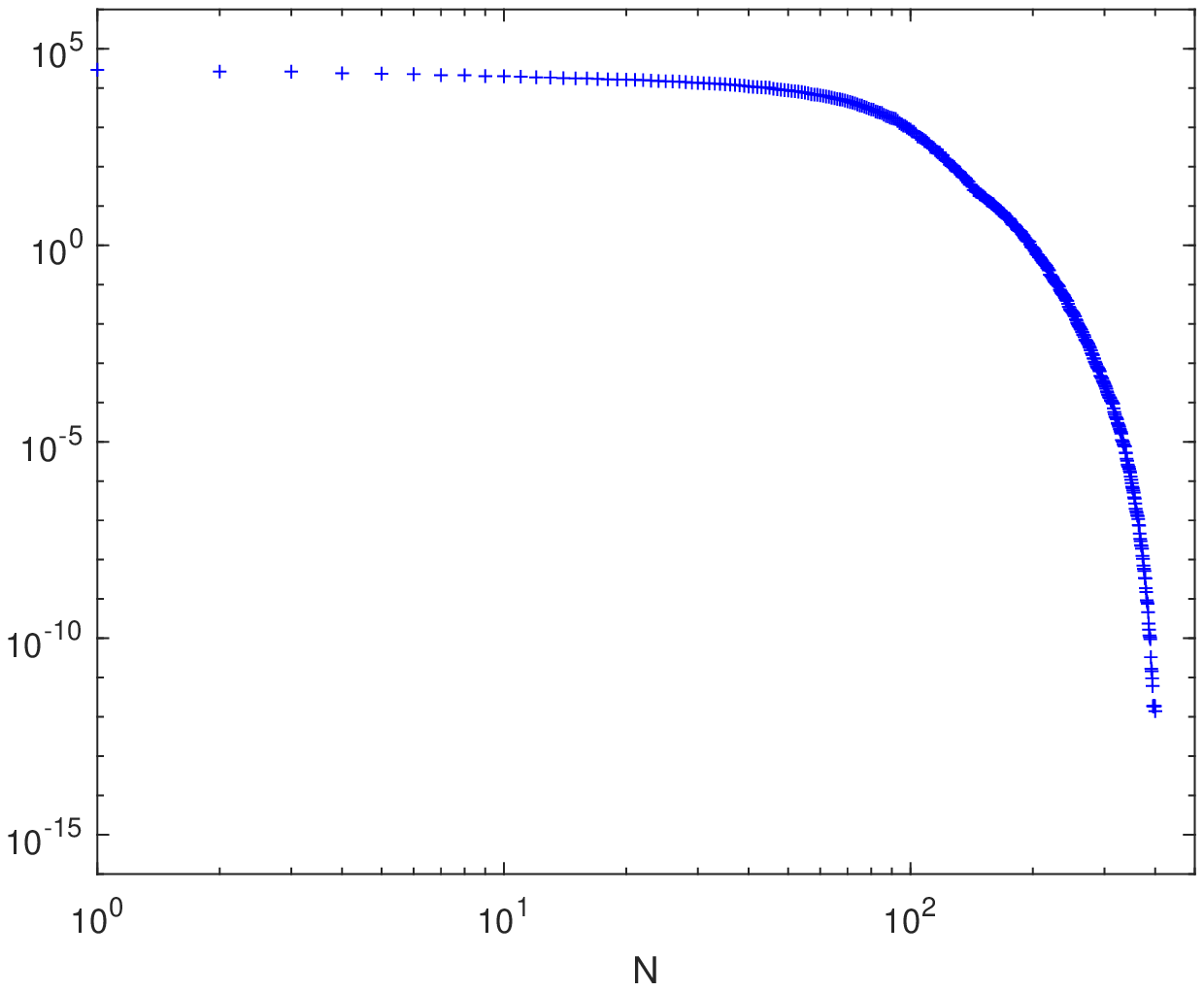}}
{\includegraphics[width=5.8cm]{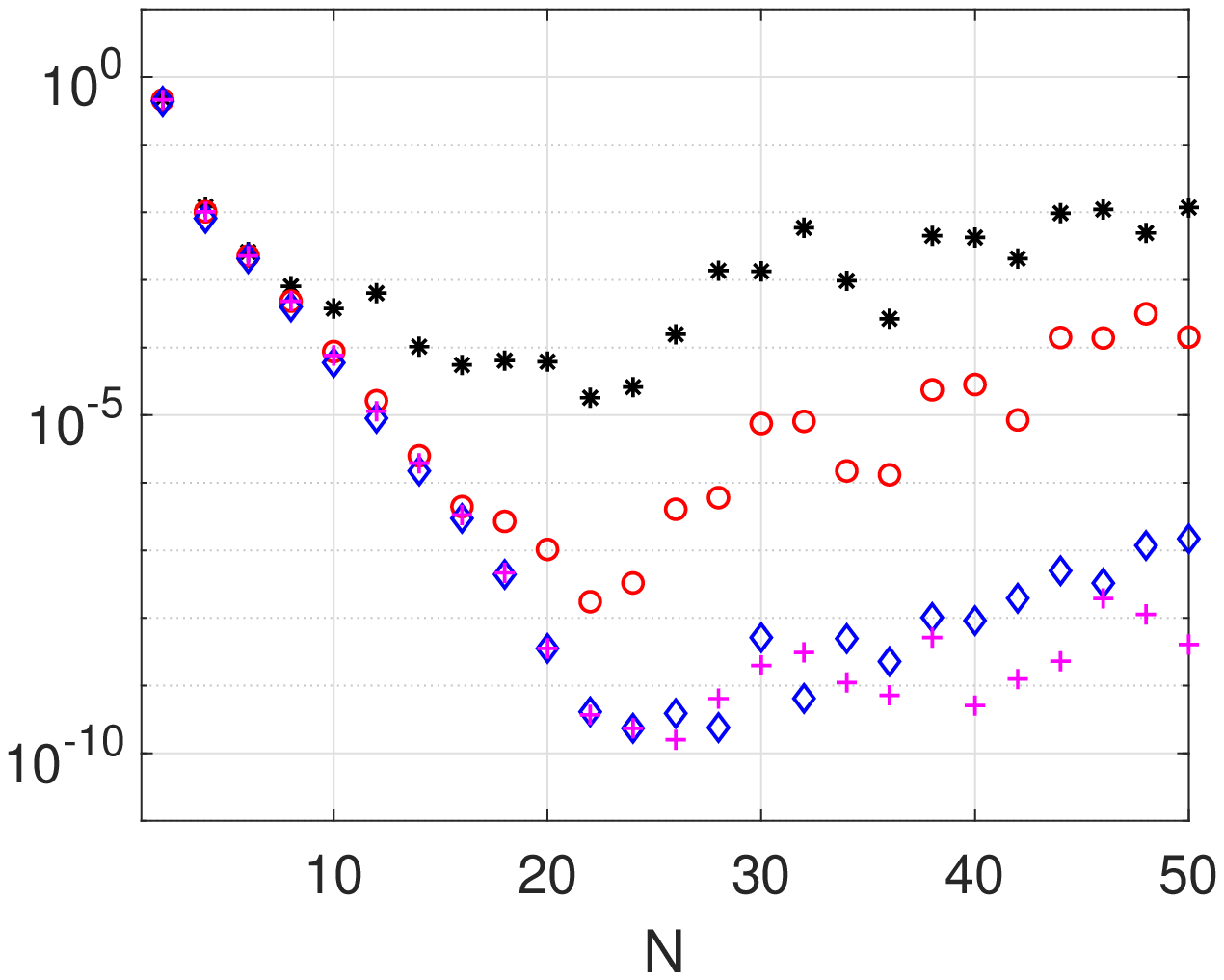}}
\caption{Left: The distribution of the singular values of $\mathbf{P}$ when $N= 20$, $N_B= 20\lfloor\log (N_\Lambda) \rfloor$. Right: The log plot of the maximum error of Example \ref{example1} with different $N_B$, i.e., $N_B= 5\lfloor\log (N_\Lambda) \rfloor$ (black star), $N_B=10 \lfloor\log (N_\Lambda) \rfloor$ (red circle), $N_B=20 \lfloor\log (N_\Lambda) \rfloor$ (blue diamond) and $N_B=5N$ (magenta plus)}
\label{PDE-EX1}
\end{figure}

In the Figure \ref{PDE-EX1}, we show the singular value profile of matrix $\mathbf{P}$ and the maximum error of PDEs with different values of $N_B$, and we omit discretization of domain $\Omega_P$. We observe that $N_B= 20 \lfloor\log (N_\Lambda) \rfloor$ nodes on the boundary are enough from Figure \ref{PDE-EX1}, and we can reduce up to hundreds of boundary nodes when compared with $N_B = 5N$. Meanwhile, we have given the approximation error of a analytic function $f_4$ on $\Omega_P$. We find that the error decay behavior of Example \ref{example1} is basically consistent with the approximation error of function $f_4$, that is, the error reaches $\mathcal{O}(10^{-10})$ and then has a slightly divergence trend.
\end{example}

\begin{example}\label{example2}
(Variable coefficient PDEs)~We set $\alpha(x,y) = \exp(x+y), \beta(x,y)= 0$  in \eqref{PDE0}. Let the exact solution be $U(x,y) = \sin\left(\pi/2\left({x^2}/{0.6^2}+{y^2}/{0.9^2}-1\right)\right)$
inside domain $\Omega_E = \left\{(x,y): {x^2}/{0.6^2}+{y^2}/{0.9^2}-1\leq0 \right\}$.

\begin{figure}[htbp]
\centering
{\includegraphics[width=5.8cm]{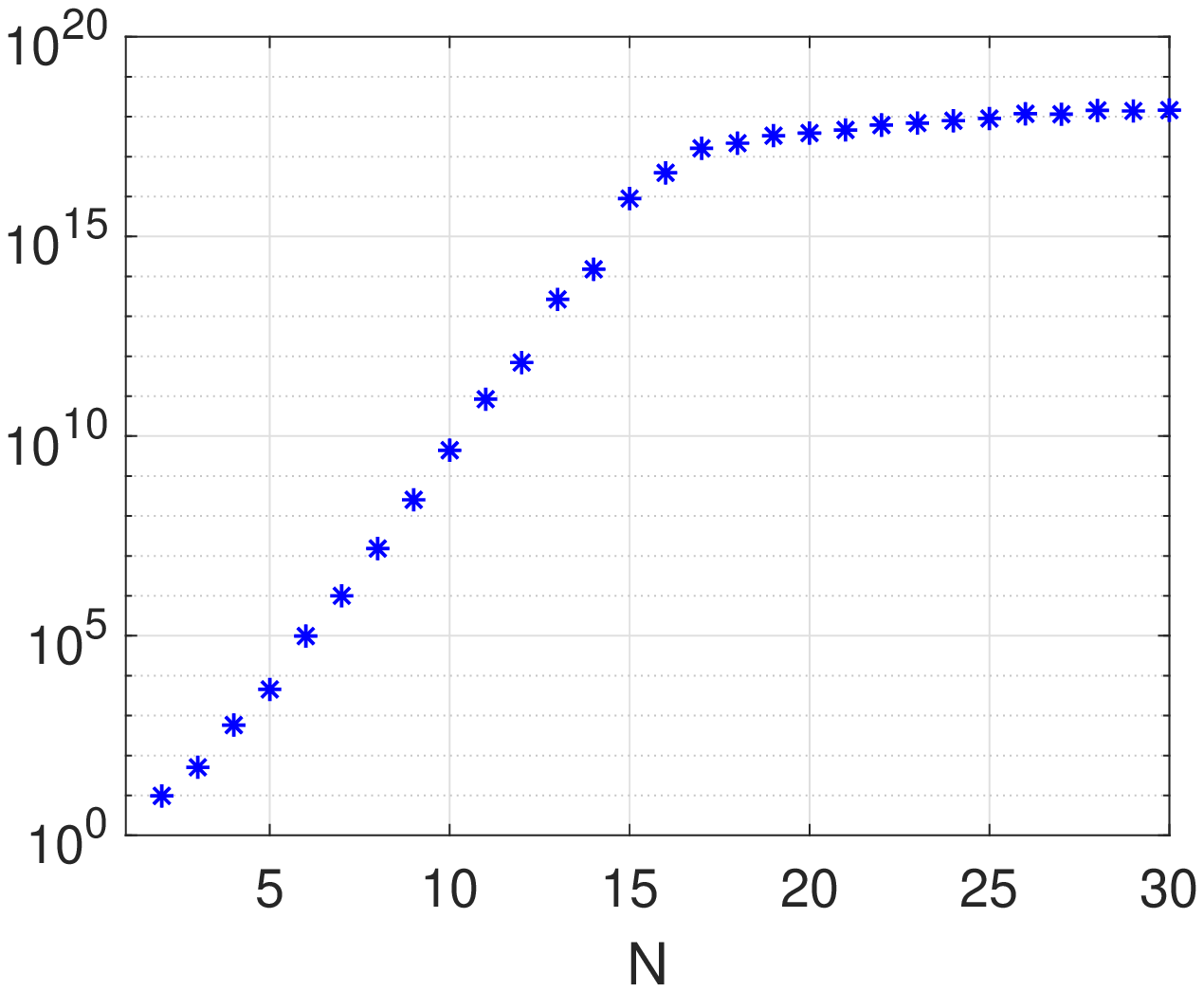}}
{\includegraphics[width=5.8cm]{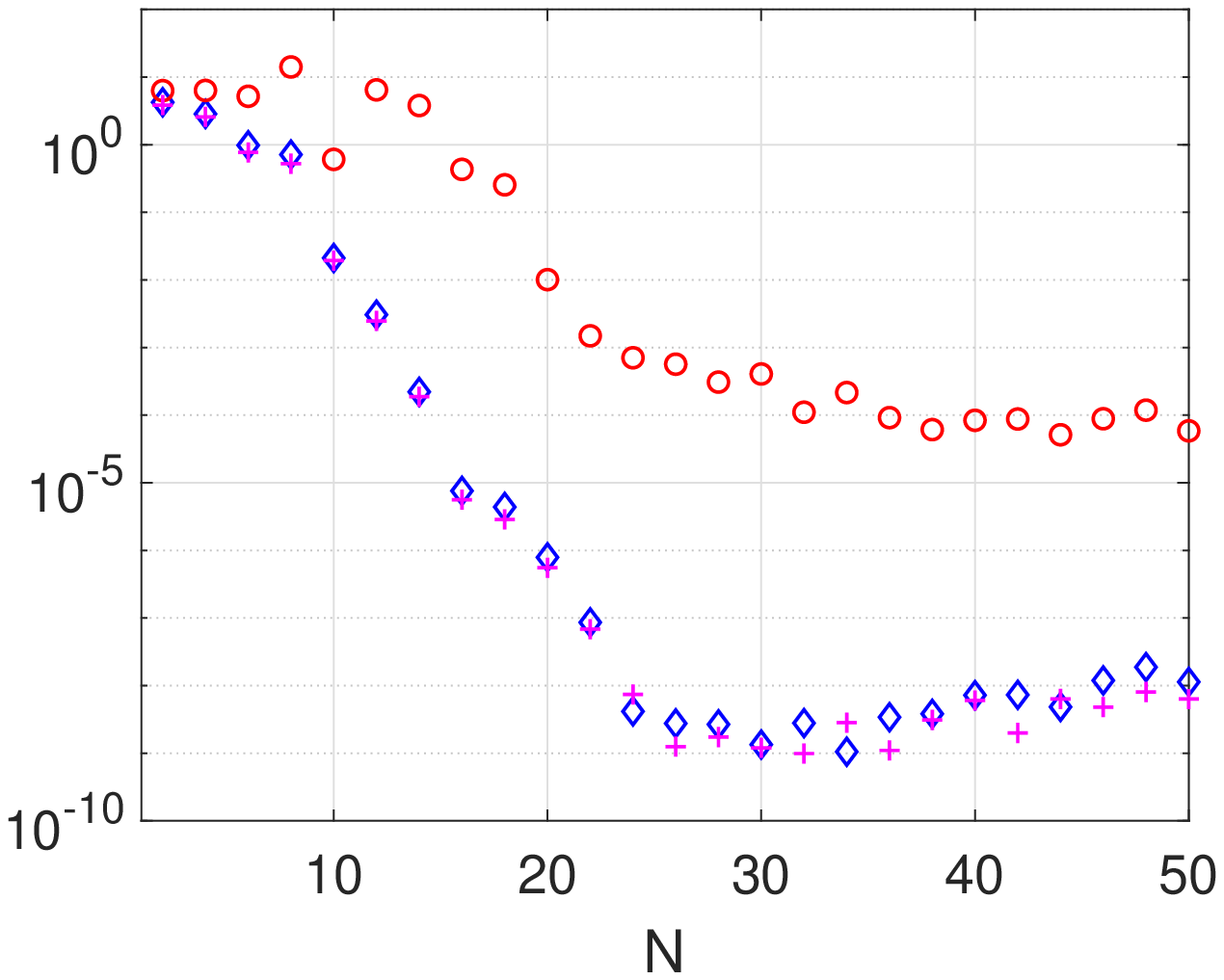}}
\caption{Left: The condition number of matrix $\mathbf{P}$ when $N= 20$, $N_B = 3N$.
Right: The log plot of the maximum error of  Example \ref{example2} with different $N_B$, i.e., $N_B=N$ (red circle), $N_B=3N$ (blue diamond) and $N_B = 5N$ (magenta plus)}
\label{PDE-EX2}
\end{figure}

The left side of Figure \ref{PDE-EX2} shows the condition number of matrix $\mathbf{P}$. This is an ill-conditioned system, and we find that the discrete systems \eqref{PDE6} of other examples also show similar ill-conditioned behavior, we will not repeat it. On right side of Figure \ref{PDE-EX2}, it shows the maximum error of Example \ref{example2} with various values of $N_B$. We observe that it is sufficient to select $N_ B=3N$, more boundary nodes do not improve the approximation accuracy at all.
\end{example}

\begin{example}\label{example3}
(Variable coefficient PDEs)~We set $\alpha(x,y) = (\sin x +1)(\cos y+1) $, $\beta(x,y) = \exp(x+y)$ in \eqref{PDE0}. Let the exact solution be $U(x,y) = \exp{\left(-(x^2+y^2)/2\right)}$ inside the triangle domain $\Omega_T$, where $\Omega_T$ has been defined in Section \ref{section2}.

\begin{figure}[htbp]
\centering
{\includegraphics[width=5.8cm]{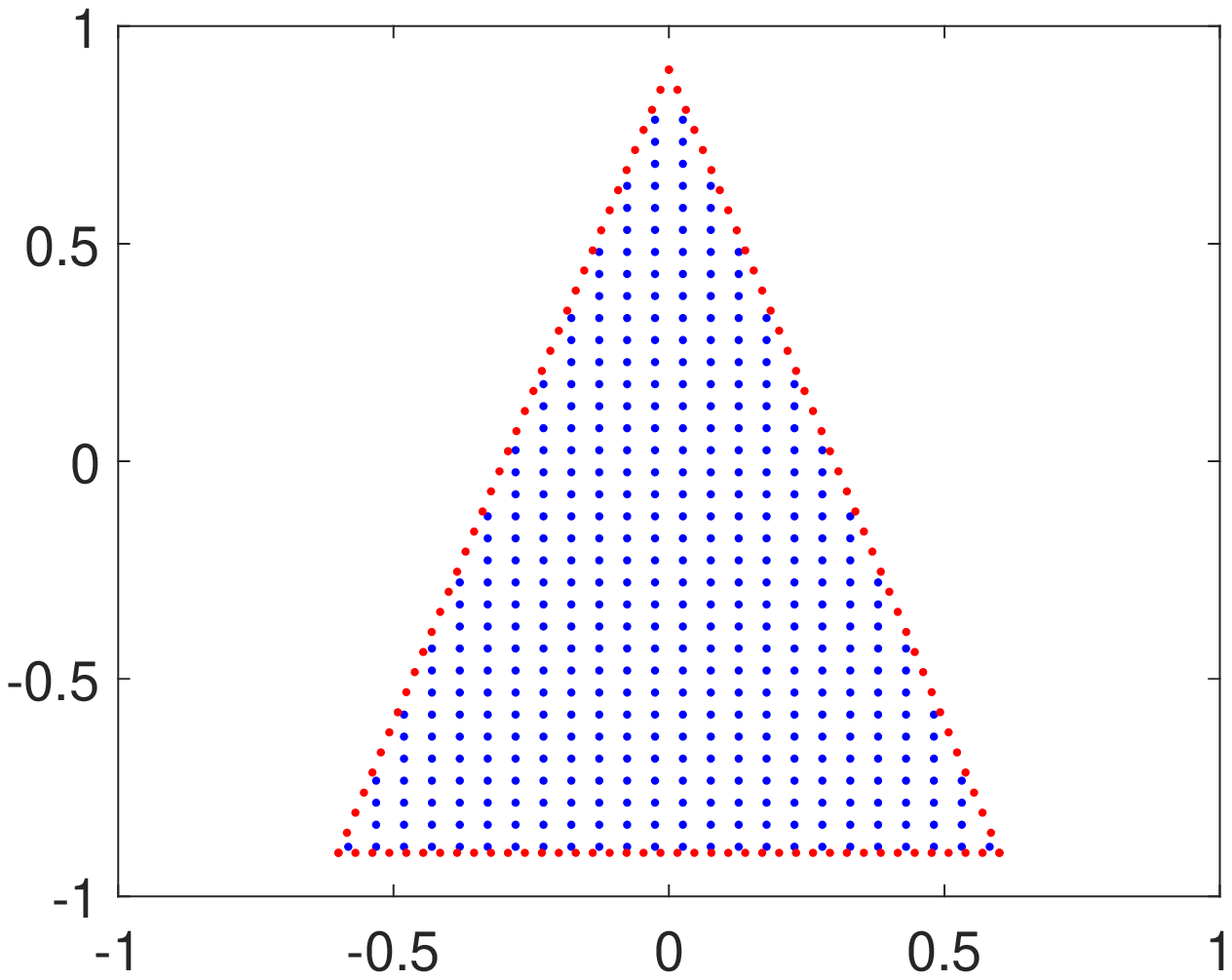}}
{\includegraphics[width=5.8cm]{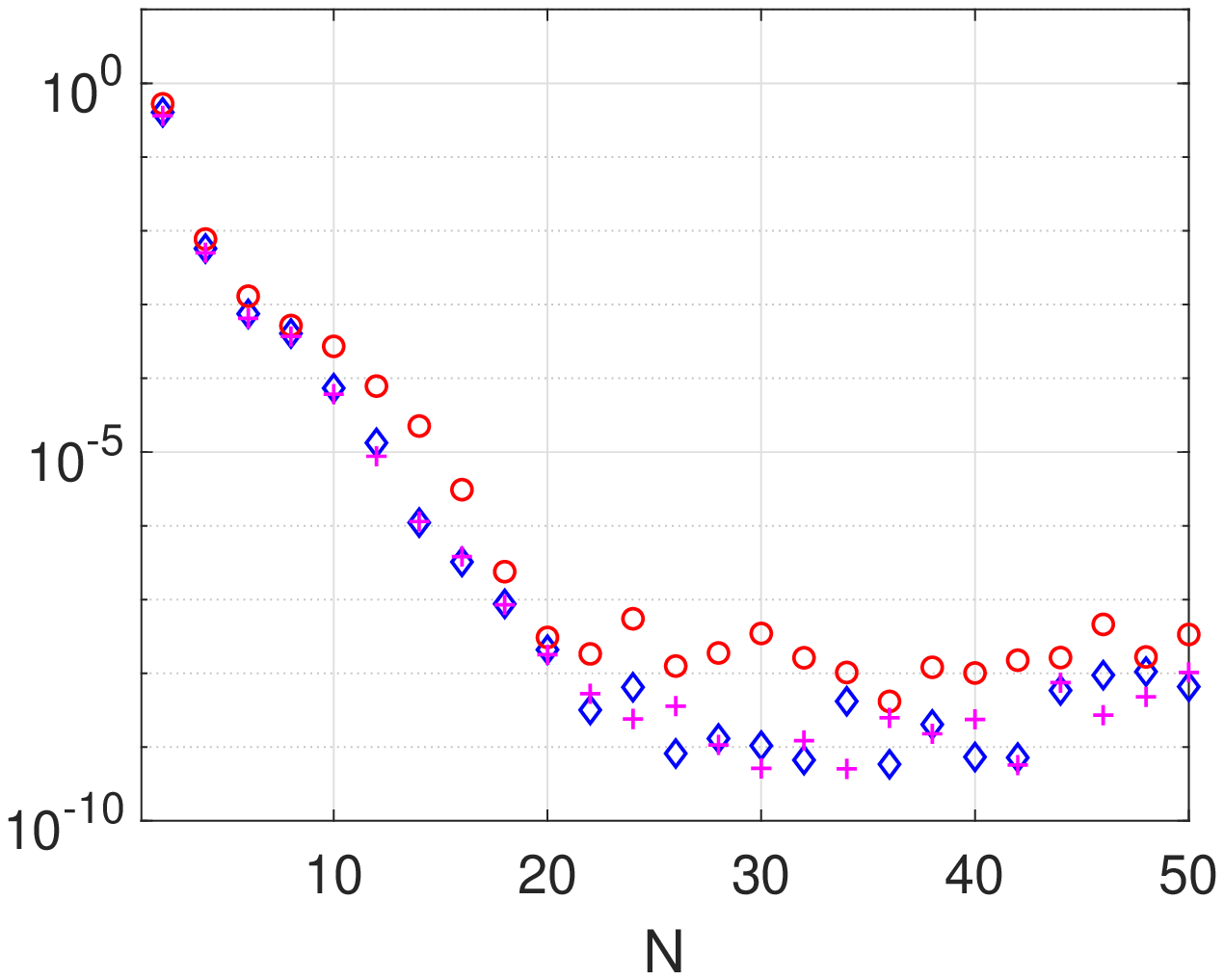}}
\caption{Left: The internal and boundary collocation nodes when $N=20$, $N_B= 6N$. Right: The log plot of the maximum error of  Example \ref{example3} with different $N_B$, i.e., $N_B=3N$ (red circle), $N_B=6N$ (blue diamond) and $N_B = 9N$ (magenta plus)}
\label{PDE-EX3}
\end{figure}
\end{example}

In Figure \ref{PDE-EX3}, we present the discretization of $\Omega_T$ when $N=20$, $N_B= 6N$, and we also show the maximum error of Example \ref{example3} with different values of $N_B$. We observe that there is no significant difference in the accuracy of these three cases, and we prefer to take $N_B = 6N$ here. We also observe that the approximation accuracy reaches about $\mathcal{O}(10^{-9})$. In Example \ref{example1}-\ref{example3}, the true solutions of these PDEs are analytic, we observe that the numerical solutions converge exponentially to a plateau, as a function of $N$. After a breakpoint, the convergence rate slows down, and there are obvious fluctuations. For analytic functions, the position of the breakpoint is almost the same, and the shape of the domain at this time does not seem to have much influences on the accuracy.

\begin{example}\label{example4}
(Corner singularity solution)~We set $\alpha(x,y) = 1, \beta(x,y) = 0$ in \eqref{PDE0}. Let the exact solution be $U(x,y) = (1-x^2)^{5/2}(1-y^2)^{5/2}$ inside a square domain $\Omega_S=\{(x,y): |x| \leq 1, |y| \leq 1\}$.

\begin{figure}[htbp]
\centering
{\includegraphics[width=5.8cm]{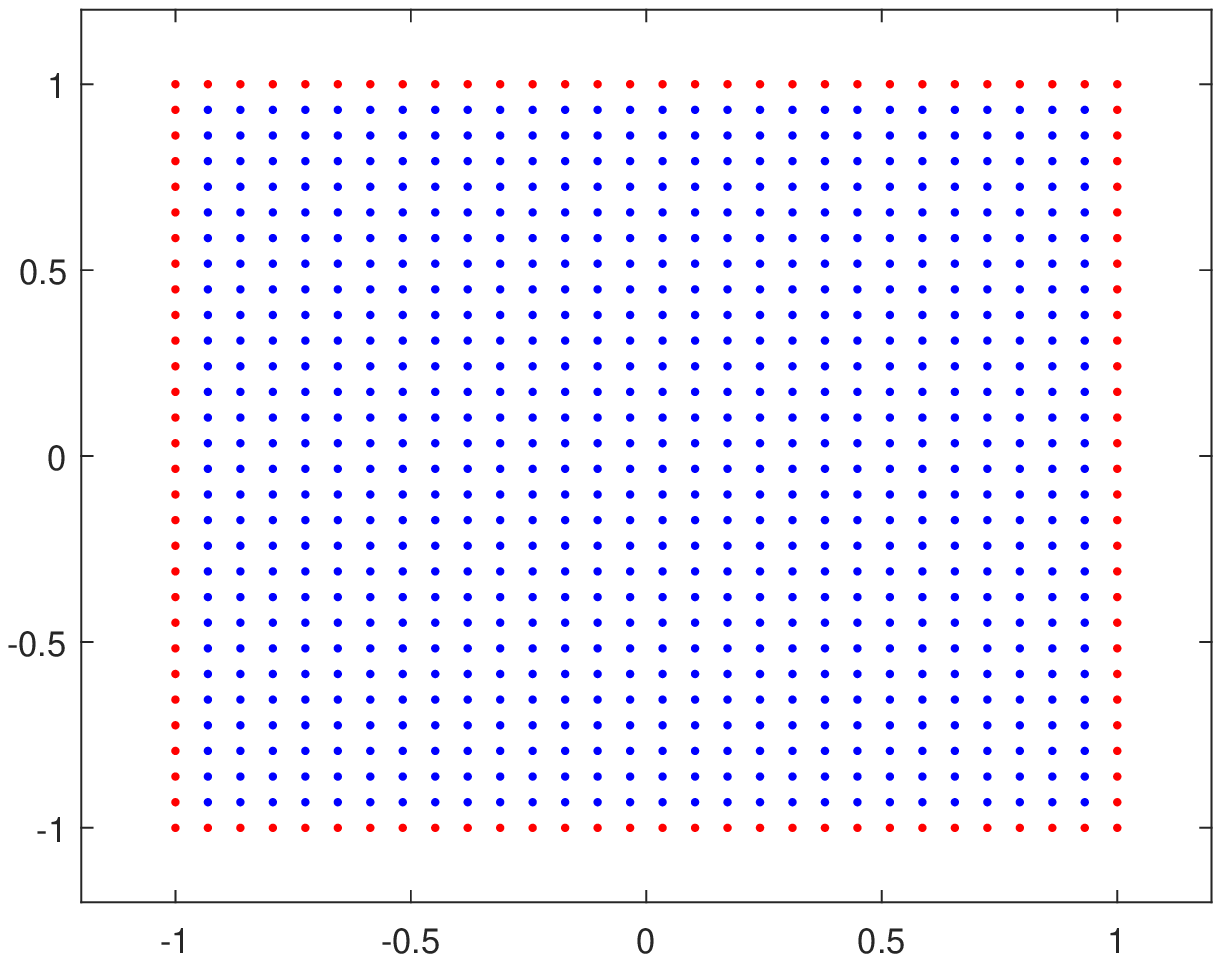}}
{\includegraphics[width=5.8cm]{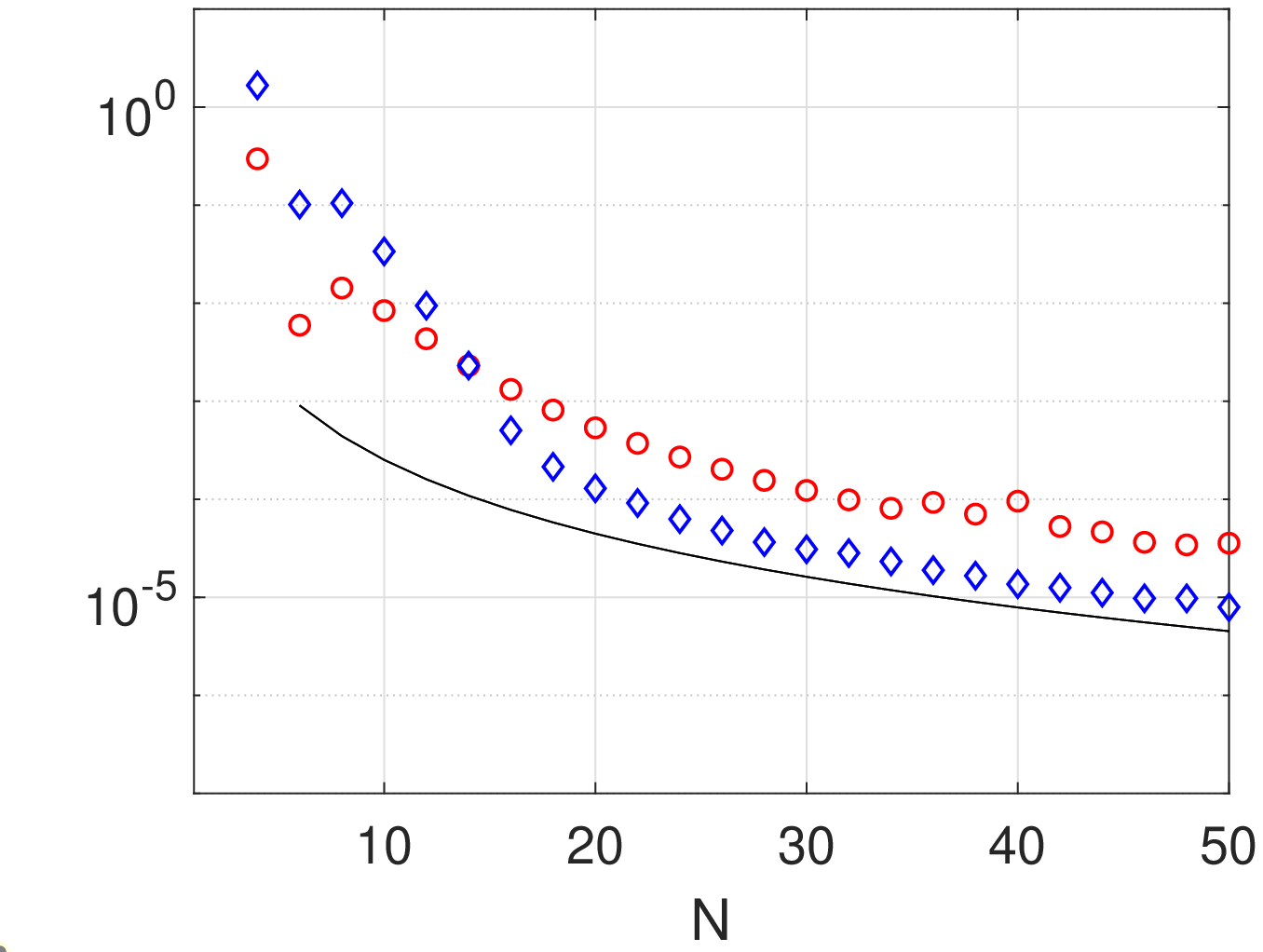}}
\caption{Left: The  internal and boundary collocation nodes when $N=30$, $N_B= 4N-4$. Right: The log plot of the maximum error of  Example \ref{example4} with different $N_B$, i.e., $N_B=4N-4$ (red circle), $N_B=8N-4$ (blue diamond). Also shown is the black curve $\mathcal{O}(N^{-5/2})$}
\label{PDE-EX4}
\end{figure}

For this special domain $\Omega_S$, we can directly define the tensor orthogonal polynomials to approximate the solutions. However, this is a corner singularity solution, we still consider to use the Fourier frames to deal with it, and the fast algorithm can also be applied in this constant coefficient PDE with suitable boundary nodes. On the left side of Figure \ref{PDE-EX4}, we give the discretization of $\Omega_T$ when $N=30$ and $N_B= 4N-4$. On the right side of Figure \ref{PDE-EX4}, we show the maximum error of PDEs, and we observe that the error graph and the black curve $\mathcal{O}(N^{-5/2})$ remain parallel when the degree of freedom is large enough.
\end{example}

\begin{example}\label{example5}
(Double connected domain)~We set $\alpha(x,y) = \exp(x+y)$, $\beta(x,y) = 0$ in \eqref{PDE0}. Let the exact solution be $U(x,y) = \sin\left(\pi(x^2+y^2-0.9^2)/2\right)$ inside domain $\Omega_{L}$, where $\Omega_{L}$ has been defined in Section \ref{section2}.

\begin{figure}[htbp]
\centering
{\includegraphics[width=5.8cm]{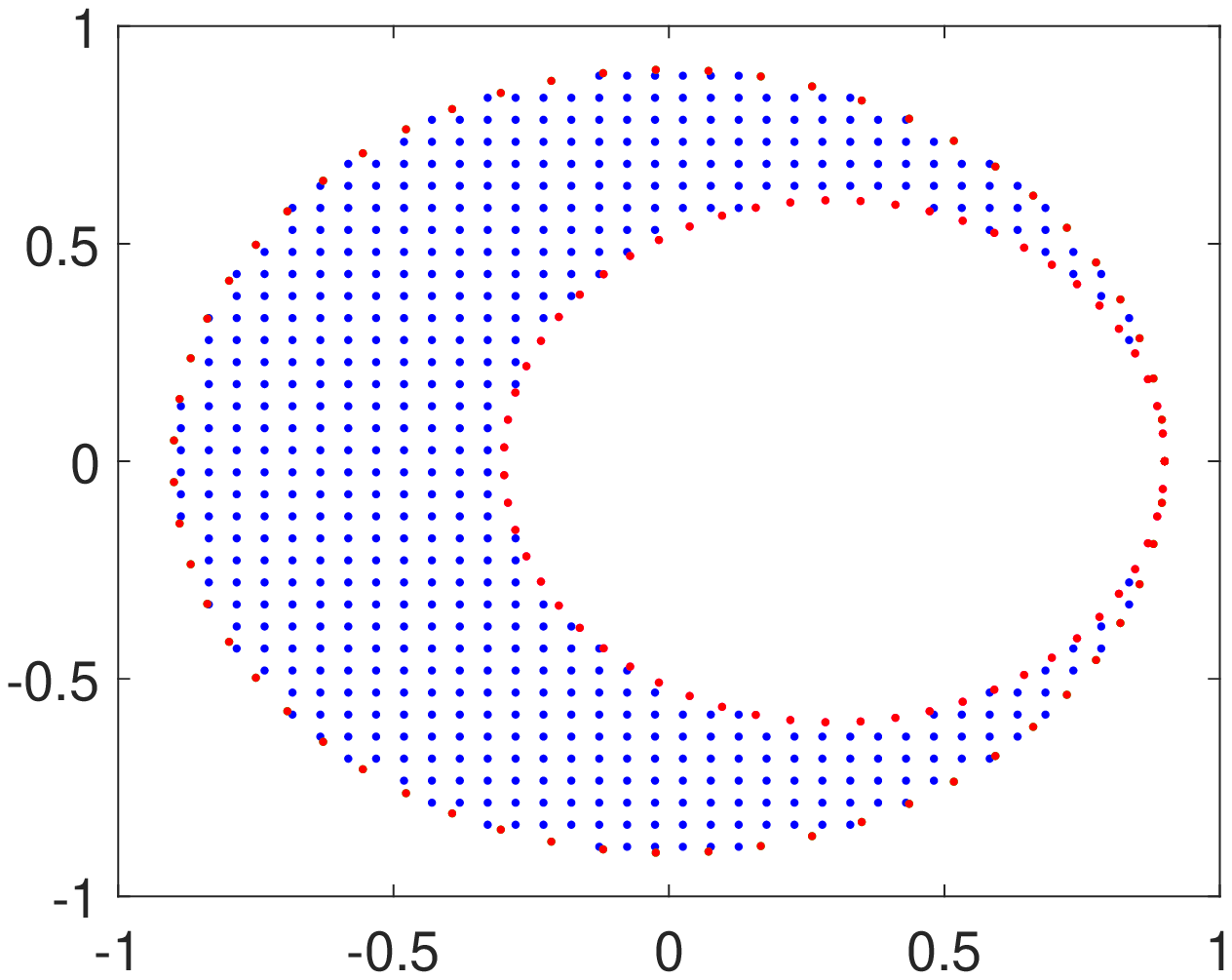}}
{\includegraphics[width=5.8cm]{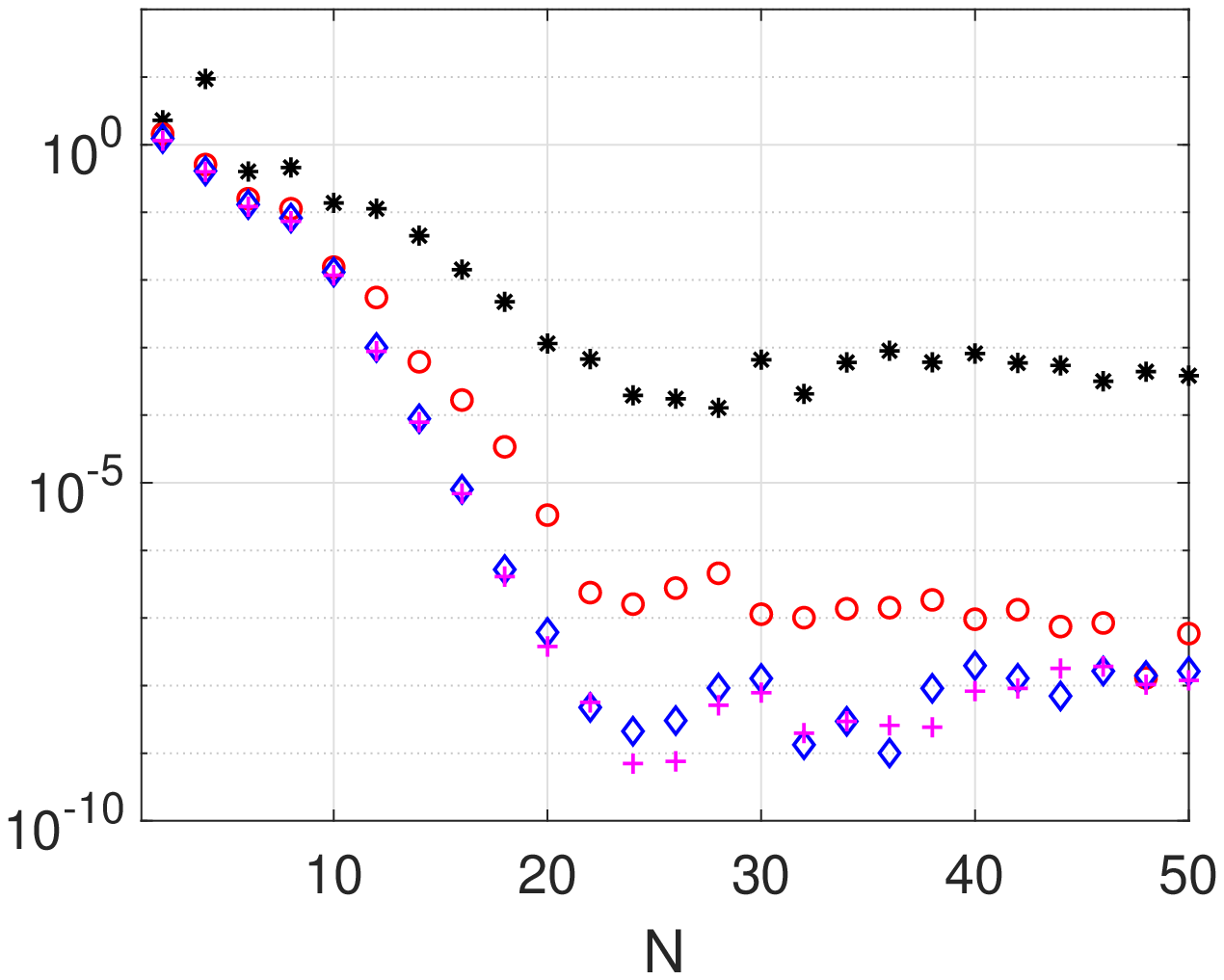}}
\caption{Left: The internal and boundary collocation nodes when $N= 20$, $N_B= 6N$. Right: The log plot of the maximum error of  Example \ref{example5}, i.e., $N_B= 2N$ (black star), $N_B= 4N$ (red circle), $N_B= 6N$ (blue diamond) and $N_B = 8N$ (magenta plus)}
\label{PDE-EX5}
\end{figure}

In Figure \ref{PDE-EX5}, we give the discretization of $\Omega_T$ when $N=20$, $N_B= 6N$, and we also show the maximum error of PDEs with $N$. The differential operator in Example \ref{example5} is the same as that in Example \ref{example2}, but is defined in different domains. One is a simply connected domain $\Omega_E$, and the other one is a doubly connected doamin $\Omega_L$. We observe that the approximation behavior of the two PDEs is almost the same, and we speculate that the connectivity of the domain will not affect the approximation accuracy of PDEs. In order to investigate the influence of domain connectivity on the accuracy of the oversampling collocation method, we investigate another doubly connected region without changing the differential operator in Example \ref{example2} and \ref{example5}.
\end{example}

\begin{example}\label{example6}
(Double connected domain)~We take $\alpha(x,y) = \exp(x+y)$, $\beta(x,y)=0$ in \eqref{PDE0}. Let the exact solution be $U(x,y) = \sin(x^2+y^2)$ inside domain $\Omega_{TF}=\{(r,\theta): 0.4+0.2\sin(5\theta) \leq r \leq  0.7+0.2\sin(5\theta)\}.$

\begin{figure}[htbp]
\centering
{\includegraphics[width=5.8cm]{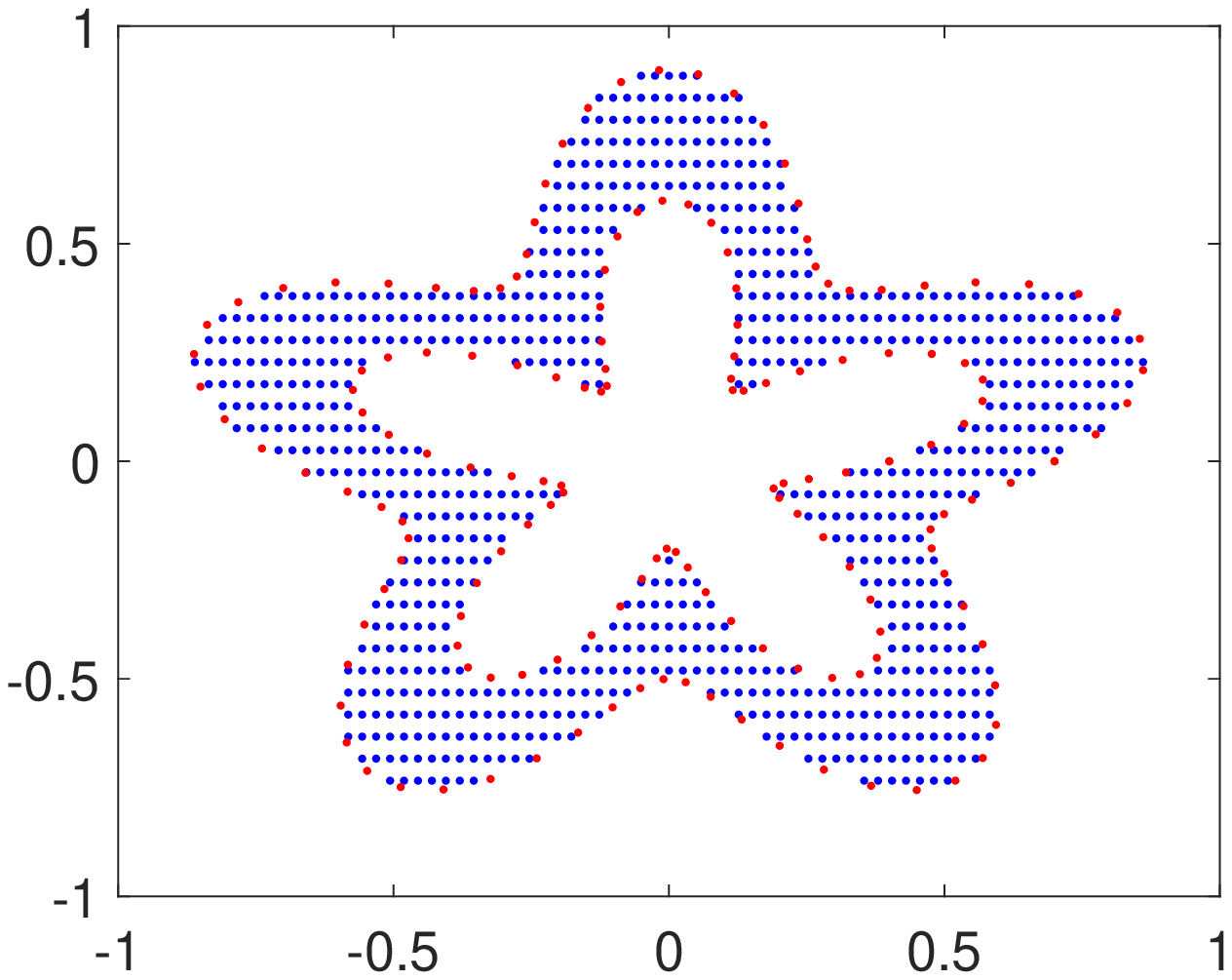}}
{\includegraphics[width=5.8cm]{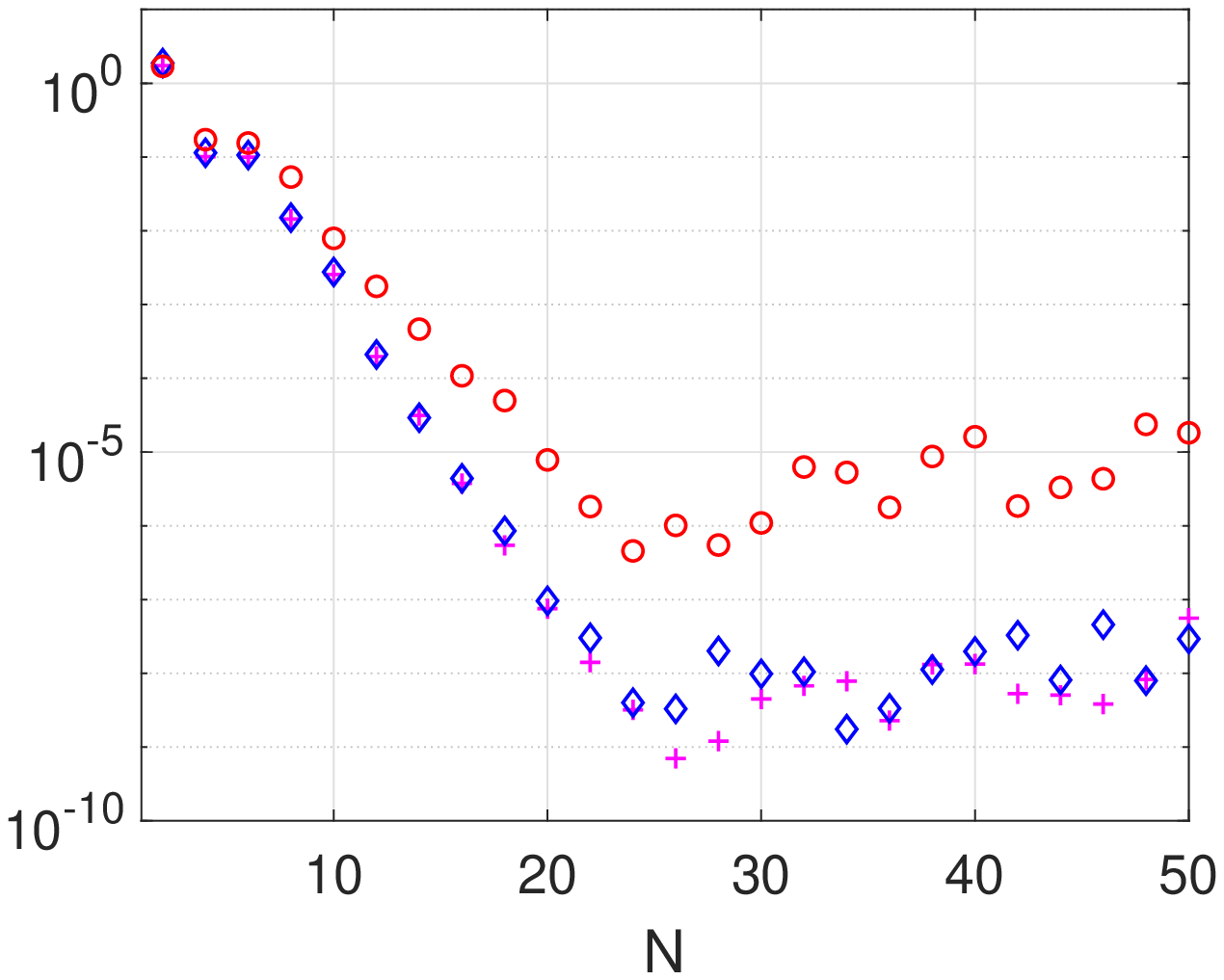}}
\caption{Left: The internal and boundary collocation nodes when $N= 20$, $N_B= 8N$. Right: The log plot of the maximum error of  Example \ref{example6}, i.e., $N_B= 4N$ (red circle), $N_B= 8N$ (blue diamond), $N_B= 10N$ (magenta plus)}
\label{PDE-EX6}
\end{figure}

On the left side of Figure \ref{PDE-EX6}, we give the discretization of $\Omega_{TF}$ when $N=20$ and $N_B= 8N$. Here, the division in the $x$-axis direction is twice as dense as that in the y-axis. On the right side of Figure \ref{PDE-EX6}, we show the maximum error with $N$. For this kind of hollowed out double connected region, we need more boundary node information to ensure accuracy, it is better to take $N_B= 8N$. Observing Figure \ref{PDE-EX5} and \ref{PDE-EX6}, we conclude that the accuracy of the collocation method is not affected, even if the region is doubly connected. The accuracy fluctuates back and forth between $\mathcal{O}(10^{-9})$ and $\mathcal{O}(10^{-8})$.
\end{example}

\begin{example}\label{example7}
(Random nodes)~We set $\alpha(x,y) = \exp(x+y)$, $\beta(x,y) = 0$ in \eqref{PDE0}. Let the exact solution be $U(x,y) = \sin\left( 4x^2- 4x^4/0.9^2  -y^2\right)$ inside
$\Omega_B = \{(x,y): 4x^2-4x^4/0.9^2 -y^2 \geq 0 \}$.

\begin{figure}[htbp]
\centering
{\includegraphics[width=5.8cm]{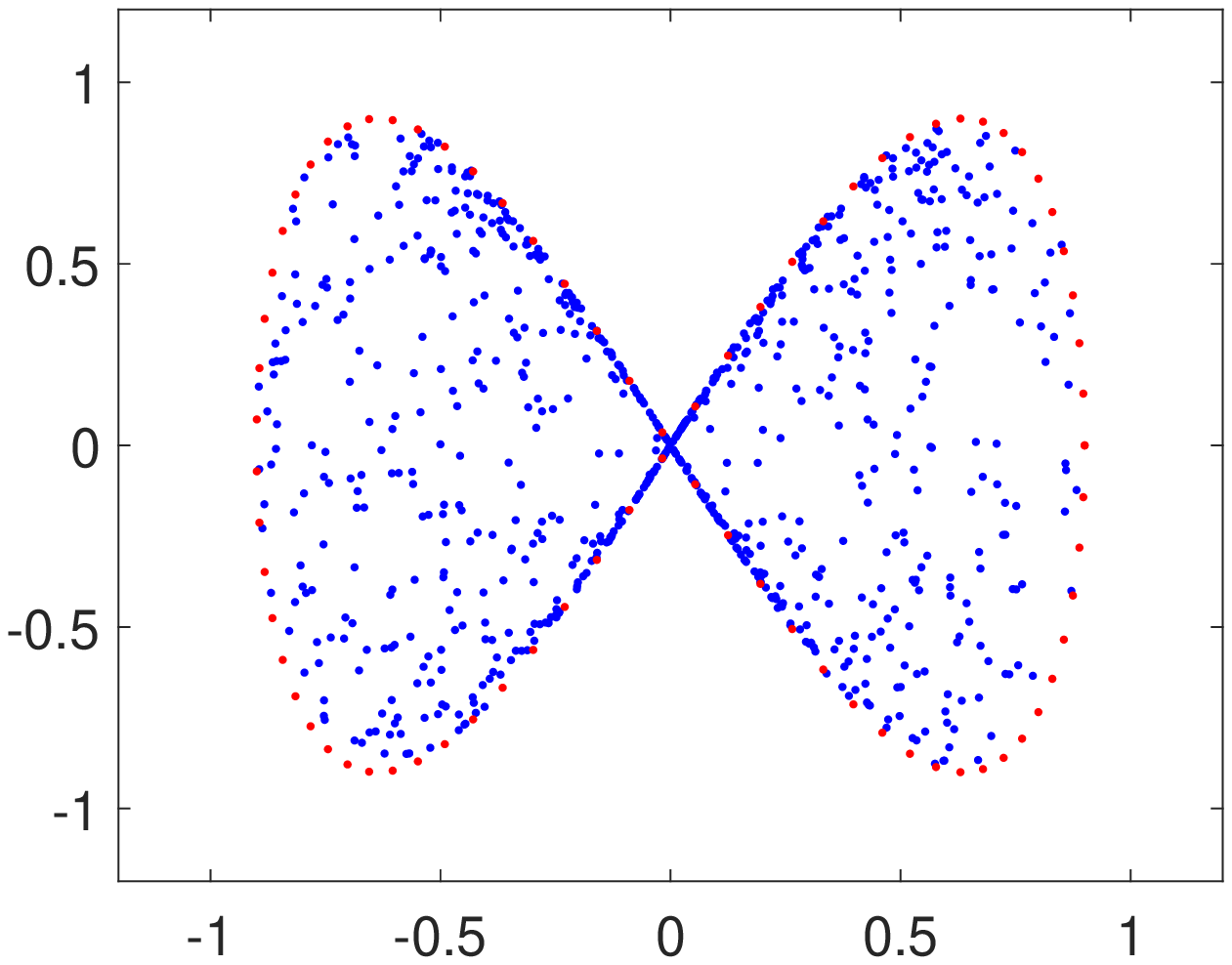}}
{\includegraphics[width=5.8cm]{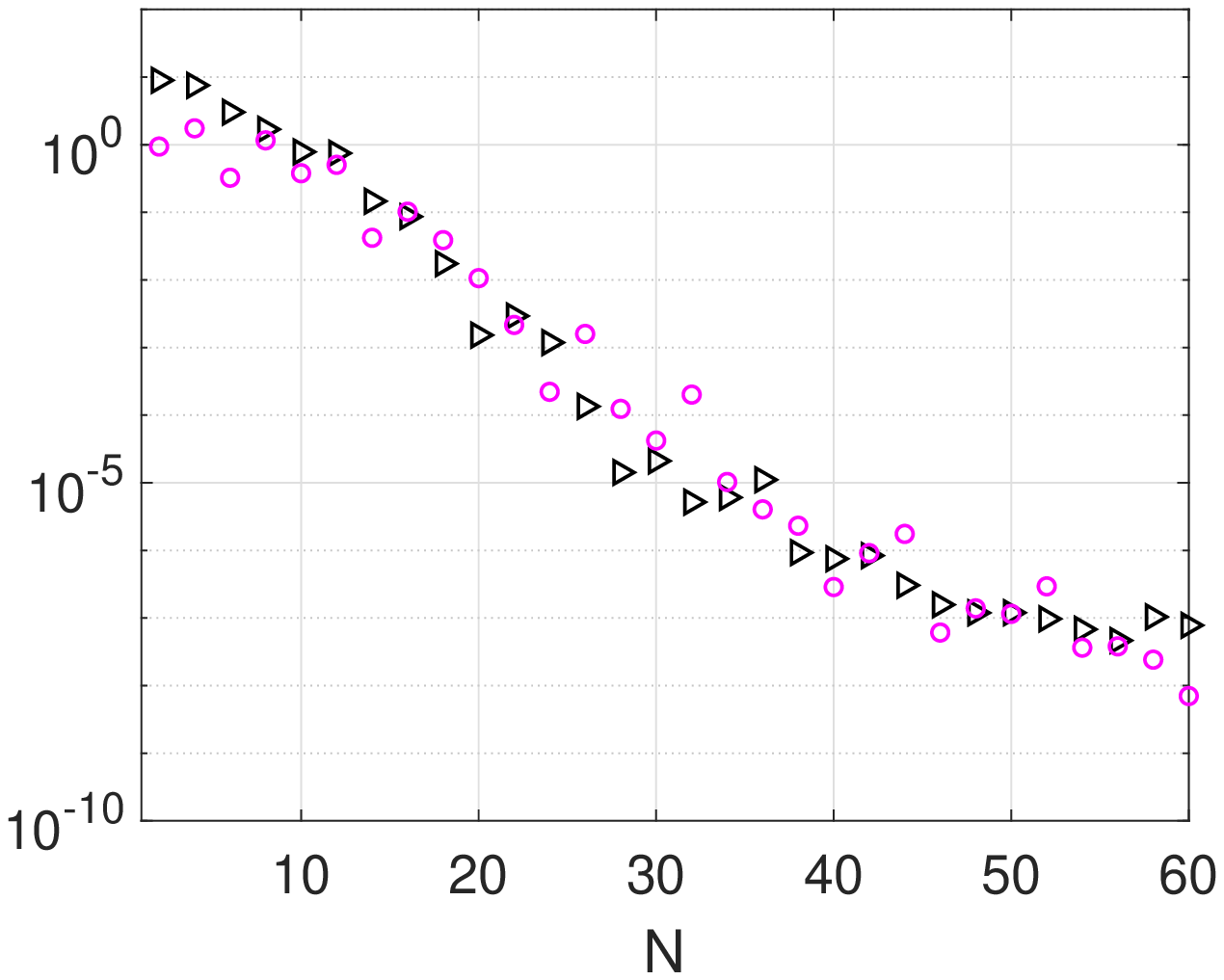}}
\caption{Left: The random internal and boundary collocation nodes when $N=20$, $N_ B=4N$. Right: The log plot of the maximum error of function $U(x,y)$ (magenta circle), and the log plot of the maximum error of Example \ref{example7} (black triangle)}
\label{PDE-EX7}
\end{figure}

Since the number of equispaced nodes strictly depends on the shape of irregular domains and division criterions, we can directly generate random nodes which satisfy the uniform distribution, then double oversampling can be achieved. On the left side of Figure \ref{PDE-EX7}, we give the discretization of $\Omega_{B}$,  here we take $N_I = 2N_\Lambda$, $N_ B=4N$. On the right side of Figure \ref{PDE-EX7}, we show the function approximation error of $U(x,y)$ and the PDE numerical approximation error. We find that the accuracy obtained by using uniform random nodes to approximate the real solution $U(x,y)$ and to solve the PDE show the same decay behavior. These two approximations have not reached equilibrium until $N=60$, although this is an analytic solution.
\end{example}

\section{Conclusions}\label{section5}
In this paper, we demonstrate a spectral collocation method for general second-order elliptic PDEs by using Fourier frames, a large number of numerical experiments show that our proposed numerical method is straightforward and performs well. But what kind of nodes to choose, how many nodes are optimal, these issues are worth to be researched. For the analytical solutions, the error decays exponentially to about $\mathcal{O}(10^{-9})$ until it reaches a breakpoint. After this point, the errors show a slightly fluctuation behavior. When we focus on the variable coefficients PDEs, the coefficients terms will change the singular value profile of the collocation matrix, then the fast algorithm of Fourier extension cannot be extended.

\section*{Declarations}
\begin{itemize}
\item {\bf Funding}: No funding was received to assist with the preparation of this manuscript.
\item {\bf Competing interests}: The authors declare no conflict of interest.
\item {\bf Authors' contributions }: These authors contributed equally to this work.
\end{itemize}

\end{document}